\documentclass[reqno, 11pt, a4paper]{amsart}

\usepackage[utf8]{inputenc}
\usepackage{amsthm, amssymb, amsmath, amsfonts}
\usepackage{url}
\usepackage[dvips]{graphicx} 
\usepackage{psfrag}
\usepackage{esint}
\usepackage{color}
\usepackage{mathrsfs}

\usepackage{array}

\usepackage[body={15cm,21.5cm},centering]{geometry} 

\newtheorem{thm}{Theorem}[section]

\newtheorem{lem}[thm]{Lemma}

\newtheorem{defi}[thm]{Definition}
\theoremstyle{remark}
\newtheorem{rem}[thm]{Remark}

\newcommand{\mb}{\mathbb}
\renewcommand{\le}{\leqslant}
\renewcommand{\ge}{\geqslant}

\newcommand{\mcl}{\mathcal}
\newcommand{\E}{\mathbb{E}}
\newcommand{\B}{\mathbb{B}}

\newcommand{\EEo}{\mathbf{E}^\omega}

\newcommand{\ee}{\mathbf{e}}
\newcommand{\N}{\mathbb{N}}
\renewcommand{\L}{\mathcal{L}}

\newcommand{\1}{\mathbf{1}}
\newcommand{\R}{\mathbb{R}}

\newcommand{\Z}{\mathbb{Z}}
\renewcommand{\P}{\mathbb{P}}
\newcommand{\PP}{\mathbf{P}}
\newcommand{\PPo}{\mathbf{P}^\omega}
\newcommand{\ov}{\overline}
\newcommand{\td}{\tilde}
\newcommand{\eps}{\varepsilon}
\def\d{{\mathrm{d}}}

\newcommand{\Var}[1]{\mathbb{V}\mathrm{ar}\left[#1\right]}
\newcommand{\mfk}{\mathfrak}

\newcommand{\Ah}{A_\mathrm{hom}}
\newcommand{\Ahd}{A_\mathrm{hom}^\mathrm{disc}}
\newcommand{\om}{\omega}
\newcommand{\diag}{\mathrm{diag}}
\newcommand{\spa}{\mathrm{spa}}
\newcommand{\law}{\mathrm{law}}

\newcommand{\ho}{\mathrm{hom}}
\newcommand{\e}{\varepsilon}
\newcommand\calF{\mathcal{F}}
\newcommand\calL{\mathcal{L}}
\newcommand{\DD}{\operatorname{D}}

\newcommand\Id{\mathrm{Id}}

\newcommand{\var}[1]{\mathrm{var}\left[#1\right]}
\newcommand{\expec}[1]{\left\langle #1 \right\rangle}

\newcommand{\dps}{\displaystyle}

\renewcommand{\leq}{\le}
\renewcommand{\geq}{\ge}

\newcommand{\Ll}{\left}
\newcommand{\Rr}{\right}
\renewcommand{\log}{\ln}

\newcommand{\osc}[2]{\underset{\dps #1}{\mathrm{osc}} \,#2\,}

\title[Convergence rates in stochastic homogenization]{
Random walk in random environment, corrector equation, and homogenized coefficients: from theory to numerics, back and forth}
\author{A.-C. Egloffe, A.~Gloria, J.-C. Mourrat, and T. N. Nguyen}

\address{A.-C. Egloffe\\ Project-team REO \\ INRIA Paris-Rocquencourt \\ Le Chesnay, France}
\email{anne-claire.egloffe@inria.fr}
\address{A. Gloria \\ Universit\'e Libre de Bruxelles (ULB) \\ Brussels, Belgium \\ and  Project-team SIMPAF \\  INRIA Lille - Nord Europe \\ Villeneuve d'Ascq, France}
\email{agloria@ulb.ac.be}
\address{J.-C. Mourrat \\ EPFL - Institut de mathématiques \\ Lausanne, Switzerland}
\email{jean-christophe.mourrat@epfl.ch}
\address{T. N. Nguyen \\ CMAP \\ \'Ecole polytechnique \\ Palaiseau, France}
\email{nhan.nguyen@cmap.polytechnique.fr}

\begin{document}
\maketitle

\begin{abstract}
This article is concerned with numerical methods to approximate effective coefficients in stochastic homogenization of discrete linear elliptic equations, and their numerical analysis --- which has been made possible by recent contributions on quantitative stochastic homogenization theory by two of us and by Otto.
This article makes the connection between our theoretical results and computations.
We give a complete picture of the numerical methods found in the literature, compare them in terms of known (or expected) convergence rates, and study them numerically. 
Two types of methods are presented: methods based on the corrector equation, and methods based on random walks in random environments.
The numerical study confirms the sharpness of the analysis (which it completes by making precise the prefactors, next to the convergence rates), supports some of our conjectures, and calls for new theoretical developments.
\end{abstract}

\medskip
{\small
\noindent \textbf{Keywords:}  stochastic homogenization, discrete elliptic equations,
effective coefficients, random walk, random environment, Monte-Carlo method, quantitative estimates.

\smallskip

\noindent \textbf{2010 Mathematics Subject Classification:} 35B27, 39A70; 60K37, 60H25, 65C05,
60H35, 60G50, 65N99.}

\bigskip

\section{Introduction}

There have always been strong connections between stochastic homogenization of linear elliptic PDEs and
 random walks in random environments (RWRE).
In the first case, the central object is the random elliptic operator which can be replaced on large scales by a deterministic elliptic operator characterized by the so-called homogenized matrix. In the second case, the distribution of the rescaled random walk approaches that of a Brownian motion, whose covariance matrix is deterministic (and corresponds to twice the homogenized matrix).

\medskip

Both approaches allow one to devise numerical methods to approximate the homogenized coefficients. Until very recently, stochastic homogenization in general was only a \emph{qualitative} theory, so that there was no quantitative analysis of these numerical methods.
In a series of papers by two of us, Otto, and Neukamm, we developed a \emph{quantitative} theory of stochastic homogenization \cite{Gloria-Otto-09,Gloria-Otto-09b,Mourrat-10,Gloria-Neukamm-Otto-a} which allowed us to make a quantitative numerical analysis of these methods \cite{Gloria-10,Gloria-Mourrat-10a,Gloria-Mourrat-10b,Gloria-Neukamm-Otto-b}.

\medskip

In the case of stochastic homogenization of discrete elliptic equations with stationary random conductances  
of finite correlation length, and the corresponding random walk on $\Z^d$, it is now possible to give a complete quantitative picture of numerical methods to approximate the homogenized coefficients.

\medskip

The aim of this article is twofold. First we address both approaches of stochastic homogenization together: the corrector equation for the PDE approach, and the RWRE for the probabilistic approach. Since the mathematics behind both approaches are somewhat different, they are rarely presented together, whereas for the quantitative analysis both can be needed (see for instance \cite{Gloria-Mourrat-10a,Gloria-Mourrat-10b}).
A crucial role is played by spectral analysis to connect both worlds, and for the quantitative analysis.
We shall present and compare several ways of approximating homogenized coefficients, based on the corrector equation and on the RWRE. For each strategy we provide the current state of the art of the numerical analysis, and motivate conjectures when the analysis is not complete.
Second, we report a systematic study of these numerical strategies. Our numerical results confirm the sharpness of the numerical analysis, support some of our conjectures, and call for new theoretical insight.

\medskip

For completeness, and in order to make this article accessible to probabilists, analysts, and numerical analysts, we provide a short comprehensive review of the well-known qualitative results of stochastic homogenization, as well as informal arguments for the quantitative results. We believe that this article, which sets and compares on rigorous ground various numerical methods on the prototypical example of stochastic homogenization of discrete elliptic equations, will also be valuable to practitioners.

\section{Stochastic homogenization: corrector equation and RWRE}
\label{s:theo}
\setcounter{equation}{0}

We start with the description of the discrete set of diffusion coefficients, first present the discrete
elliptic point of view, and then turn to the random walk in random
environment viewpoint.
The aim of this section is to introduce a formalism, and give an
intuition on both points of view. 

\medskip

The results recalled here are essentially due to Papanicolaou and Varadhan \cite{Papanicolaou-Varadhan-79}, Kozlov \cite{Kozlov-87}, and
Kipnis and Varadhan \cite{Kipnis-Varadhan-86}).

\medskip

We present the results in the case of independent and identically distributed (i.i.d.)\ conductivities, although everything in this section remains valid provided the conductivities lie in a compact set of $(0,+\infty)$,
are stationary, and ergodic.
In particular, we shall apply this theory (and its quantitative counterpart) to coefficients which have finite correlation length. Yet, for the sake of clarity, we stick to i.i.d.\ coefficients in this presentation.

\subsection{Random environment}

We say that $x, y$ in $\Z^d$ are neighbors, and write $x \sim y$, whenever $|y-x| = 1$. This relation turns $\Z^d$ into a graph, whose set of (non-oriented) edges is denoted by $\mb{B}$.
We now define the conductances on $\mb{B}$, and their statistics.
\begin{defi}[environment]\label{2-1:defi:envi}
Let $0 < \alpha \le \beta < + \infty$, and $\Omega = [\alpha,\beta]^{\mb{B}}$. An element $\om = (\om_e)_{e \in \mb{B}}$ of $\Omega$ is called an \emph{environment}. 
With any edge $e = (x,y) \in \mb{B}$, we associate the \emph{conductance} $\om_{(x,y)}:=\om_e$ (by construction $\om_{(x,y)} = \om_{(y,x)}$). 
Let $\nu$ be a probability measure on $[\alpha,\beta]$. 
We endow $\Omega$ with the product probability measure $\P = \nu ^{\otimes \mb{B}}$. 
In other words, if $\omega$ is distributed according to the measure $\P$, then $(\omega_e)_{e \in \mb{B}}$ are independent random variables of law $\nu$. 
We denote by $L^2(\Omega)$ the set of real square integrable functions on $\Omega$ for the measure $\mathbb{P}$,
and write $\expec{\cdot}$ for the expectation associated with $\P$. 
\end{defi}
We then introduce the notion of stationarity.
\begin{defi}[stationarity]\label{2-1:defi:stat}
For all $z\in \Z^d$, we let $\theta_z:\Omega\to \Omega$ be such that for all $\om \in \Omega$ and $(x,y)\in \mb{B}$,
$(\theta_z \ \omega)_{(x,y)}=\omega_{(x+z,y+z)}$. This defines an additive action group $\{\theta_z\}_{z\in \Z^d}$ on $\Omega$ which preserves
the measure $\mathbb{P}$, and is ergodic for  $\mathbb{P}$.

We say that a function $f:\Omega \times \Z^d\to \R$ is \emph{stationary} if and only if for all $x,z\in\Z^d$ and
$\mathbb{P}$-almost every $\omega\in \Omega$, 
$$
f(x+z,\omega)\,=\,f(x,\theta_z \ \omega).
$$
In particular, with all $f\in L^2(\Omega)$, one may associate the stationary function (still denoted by $f$) : $\Z^d\times\Omega\to \R,
(x,\omega) \mapsto f(\theta_x \ \omega)$.
In what follows we will not distinguish between $f\in L^2(\Omega)$ and its stationary extension on $\Z^d\times\Omega$.
\end{defi}

\subsection{Corrector equation}\label{sec:2-1-1-2}

We associate with the conductances on $\mb B$ a conductivity matrix on $\Z^d$.
\begin{defi}[conductivity matrix]\label{2-1:defi:conduct}
Let $\Omega$, $\mathbb{P}$, and $\{\theta_z\}_{z\in \Z^d}$ be as in Definitions~\ref{2-1:defi:envi} and~\ref{2-1:defi:stat}.
The stationary diffusion matrix $A:\Z^d\times\Omega\to \mathcal{M}_d(\R)$ is defined by
$$
A:(x,\omega) \mapsto \diag(\om_{(x,x+\ee_i)} ,\dots, \om_{(x,x+\ee_d)}).
$$
\end{defi}
For each $\omega\in\Omega$, we consider the discrete elliptic operator $L$ defined by
\begin{equation}\label{2-1:eq:def-op-space}
L=-\nabla^*\cdot A (\cdot,\omega)\nabla,
\end{equation}
where $\nabla$ and $\nabla^*$ are the forward and backward discrete gradients, acting on functions $u:\Z^d\to \R$ by
\begin{equation}\label{2-1:eq:disc-nabla}
\nabla u(x):=\left[  
\begin{array}{l}
u(x+\ee_1)-u(x) \\
\vdots\\
u(x+\ee_d)-u(x)
\end{array}
\right],
\
\nabla^* u(x):=\left[  
\begin{array}{l}
u(x)-u(x-\ee_1) \\
\vdots\\
u(x)-u(x-\ee_d)
\end{array}
\right],
\end{equation}
and we denote by $\nabla^*\cdot$ the backward divergence.
In particular, for all $u:\Z^d\to \R$,
\begin{equation}\label{2-1:eq:def-op-space2}
 Lu:\Z^d\to \R,\  z\,\mapsto\, \sum_{z'\sim z}\omega_{(z,z')} (u(z)-u(z')).
\end{equation}
The standard stochastic homogenization theory for such discrete elliptic operators (see for instance \cite{Kunnemann-83}, \cite{Kozlov-87}) 
ensures that 
there exist homogeneous and deterministic coefficients $A_\ho$ such that the solution operator of the continuum
differential operator $-\nabla \cdot A_\ho\nabla$ describes $\mathbb{P}$-almost surely 
the large scale behavior of the solution operator
of the discrete differential operator $-\nabla^*\cdot A (\cdot,\omega)\nabla $.
As for the periodic case, the definition of $A_\ho$ involves the so-called correctors.
Let $\xi\in \R^d$ be a fixed direction. The corrector $\phi:\Z^d\times \Omega \to \R$ in the direction $\xi$ is the unique solution (in a sense made precise below) to
the equation
\begin{equation}\label{2-1:eq:corr-sto}
-\nabla^*\cdot A(x,\omega)(\xi+\nabla \phi(x,\omega))\,=\,0,\qquad x \in \Z^d.
\end{equation}
The following lemma gives the existence and uniqueness of this corrector $\phi$.
\begin{lem}[corrector]\label{2-1:lem:corr}
Let $\Omega$, $\mathbb{P}$, $\{\theta_z\}_{z\in \Z^d}$, and $A$ be as in Definitions~\ref{2-1:defi:envi}, \ref{2-1:defi:stat}, and~\ref{2-1:defi:conduct}.
Then, for all $\xi\in \R^d$, there exists a unique measurable function $\phi:\Z^d\times \Omega\to \R$ such that 
$\phi(0,\cdot)\equiv 0$, $\nabla \phi$ is stationary, $\expec{\nabla \phi}=0$, and $\phi$ solves \eqref{2-1:eq:corr-sto} $\mathbb{P}$-almost surely.
Moreover, the symmetric homogenized matrix $A_\ho$ is characterized by
\begin{equation}\label{2-1:eq:hom-coeff}
\xi\cdot A_\ho\xi\,=\,\expec{(\xi+\nabla\phi)\cdot A(\xi+\nabla \phi)}.
\end{equation}
\end{lem}
The standard proof of Lemma~\ref{2-1:lem:corr} 
makes use of the regularization of \eqref{2-1:eq:corr-sto} by a zero-order term $\mu>0$:
\begin{equation}\label{2-1:eq:mod-corr-sto}
\mu\phi_\mu(x,\omega)-\nabla^*\cdot A(x,\omega)(\xi+\nabla \phi_\mu(x,\omega))\,=\,0, \qquad x \in \Z^d.
\end{equation}
\begin{lem}[regularized corrector]\label{2-1:lem:mod-corr}
Let $\Omega$, $\mathbb{P}$, $\{\theta_z\}_{z\in \Z^d}$, and $A$ be as in Definitions~\ref{2-1:defi:envi}, \ref{2-1:defi:stat}, and~\ref{2-1:defi:conduct}.
Then, for all $\mu>0$ and $\xi\in \R^d$, there exists a unique stationary function $\phi_\mu \in L^2(\Omega)$ with $\expec{\phi_\mu}=0$
which solves \eqref{2-1:eq:mod-corr-sto} $\mathbb{P}$-almost surely.
\end{lem}

To prove Lemma~\ref{2-1:lem:mod-corr}, we follow \cite{Papanicolaou-Varadhan-79}, and introduce difference operators on $L^2(\Omega)$: 
for all $u\in L^2(\Omega)$, we set
\begin{equation}\label{2-1:eq:disc-nabla-sto}
\DD u(\omega):=\left[  
\begin{array}{l}
u(\theta_{\ee_1}\omega)-u(\omega) \\
\vdots\\
u(\theta_{\ee_d}\om)-u(\omega)
\end{array}
\right],
\
\DD^*u(\omega):=\left[  
\begin{array}{l}
u(\om)-u(\theta_{-\ee_1}\om) \\
\vdots\\
u(\om)-u(\theta_{-\ee_d}\om)
\end{array}
\right].
\end{equation} 
These operators play the same
roles as the finite differences $\nabla$ and $\nabla^*$ --- this time for
the variable $\omega$ (in other words, they define a difference calculus on $L^2(\Omega)$). They  allow us to define the counterpart on $L^2(\Omega)$ to the operator $L$ of \eqref{2-1:eq:def-op-space}:
\begin{defi}\label{2-1:defi:operator-sto}
Let $\Omega$, $\mathbb{P}$, $\{\theta_z\}_{z\in \Z^d}$, and $A$ be as in Definitions~\ref{2-1:defi:envi}, \ref{2-1:defi:stat}, and~\ref{2-1:defi:conduct}.
We define $\calL:L^2(\Omega)\to L^2(\Omega)$ by
\begin{eqnarray*}
\calL u(\omega)&=&-\DD^*\cdot A(\omega) \DD u(\omega) \\
&=&\sum_{z \sim 0} \om_{0,z} (u(\om) -u(\theta_z \ \om))
\end{eqnarray*}
%
where $\DD$ and $\DD^*$ are as in \eqref{2-1:eq:disc-nabla-sto}.
\end{defi}

The fundamental relation between $L$ and $\calL$ is the following
identity for stationary fields $u:\Z^d\times \Omega \to \R$: for all
$z\in \Z^d$ and almost every $\omega \in \Omega$,
\begin{equation*}
Lu(z,\omega) \,=\,\calL u(\theta_z\omega).
\end{equation*}
In particular, the regularized corrector $\phi_\mu$ is also the  unique weak
solution in $L^2(\Omega)$ to the equation
\begin{equation*}
\mu \phi_\mu(\omega)-\DD^*\cdot A(\omega)(\xi+\DD
\phi_\mu(\omega))\,=\,0, \quad \omega \in \Omega,
\end{equation*}
and its existence simply follows from the Riesz representation theorem on $L^2(\Omega)$.

The regularized corrector $\phi_\mu$ is an approximation of the corrector $\phi$ in the following sense:
\begin{lem}\label{2-1:lem:corr-mod-corr}
Let $\Omega$, $\mathbb{P}$, $\{\theta_z\}_{z\in \Z^d}$, and $A$ be as in Definitions~\ref{2-1:defi:envi}, \ref{2-1:defi:stat}, and~\ref{2-1:defi:conduct}.
For all $\mu>0$ and $\xi\in \R^d$, let $\phi$ and $\phi_\mu$ be the corrector and regularized corrector of Lemmas~\ref{2-1:lem:corr} and~\ref{2-1:lem:mod-corr}.
Then, we have
\begin{equation*}
\lim_{\mu\to 0}\expec{|\DD\phi_\mu-\DD\phi|^2}\,=\,0.
\end{equation*}
\end{lem}
From the elementary a priori estimates 
\begin{equation*}
\expec{|\nabla \phi_\mu|^2}\,=\,\expec{|\DD \phi_\mu|^2}\,\leq \, C,
\qquad \expec{\phi_\mu^2}\,\leq \, C \mu^{-1},
\end{equation*}
for some $C$ independent of $\mu$, we learn that $\DD\phi_\mu$ is bounded in $L^2(\Omega,\R^d)$ uniformly in $\mu$, so that
up to extraction it converges weakly in $L^2(\Omega,\R^d)$ to some random field $\Phi$ (which is a gradient).
This allows one to pass to the limit in the weak formulations and obtain
the existence of a field $\Phi=(\Phi_1,\dots,\Phi_d) \in L^2(\Omega,\R^d)$ such that for all
$\psi\in L^2(\Omega)$,
\begin{equation}\label{eq:corr-eq-ant1}
\expec{\DD \psi \cdot A(\xi+\Phi)}\,=\,0.
\end{equation}
Using the following weak Schwarz commutation rule
\begin{equation*}
\forall j,k \in \{1,\dots,d\}, \quad \expec{(\DD_j \psi ) \Phi_k}\,=\,\expec{(\DD_k \psi)  \Phi_j}
\end{equation*}
one may define $\phi:\Z^d\times \Omega\to \R$ such that $\nabla \phi$ is stationary, $\Phi=\nabla
\phi$, and $\phi(0,\omega)=0$ for almost every $\omega\in \Omega$. By definition this function $\phi$ is not stationary.
It is a priori not clear (and even wrong in dimension $d \le 2$) whether there exists some 
function $\psi\in L^2(\Omega)$ such that $\DD \psi=\Phi$ (this is a major difference with the periodic case).

The uniqueness of $\Phi$ is a consequence of Lemma~\ref{2-1:lem:corr-mod-corr}, which follows from the fact that $\DD \phi_\mu$ is a Cauchy sequence
in $L^2(\Omega)$. To prove this, we shall appeal to spectral theory.

The operator $\calL$ of Definition \ref{2-1:defi:operator-sto}  is bounded, self-adjoint, and non-negative on $L^2(\Omega)$. Indeed, for all $\psi,\chi \in L^2(\Omega)$, we have by direct computations
$$
\expec{(\calL \psi)^2}^{1/2} \,\leq\, 4d \sqrt{\beta} \expec{\psi^2}^{1/2}, \quad 
\expec{ (\calL \psi) \chi}\,=\,\expec{\psi (\calL \chi)}, \quad \expec{\psi \calL \psi}\geq 2d \alpha \expec{\psi^2} .
$$
Hence, $\calL$ admits a spectral decomposition in $L^2(\Omega)$.
For all $g\in L^2(\Omega)$ we denote by $e_g$ the projection of the spectral measure of $\calL$ on $g$.
This defines the following spectral calculus: for any bounded continuous function $\Psi:[0,+\infty)\to \R_+$, 
\begin{equation*}
\expec{(\Psi(\calL)g)g}\,=\,\int_{\R^+}\Psi(\lambda)de_g(\lambda). 
\end{equation*}
Let $\xi\in \R^d$ with $|\xi|=1$ be fixed, and define the local drift as $\mathfrak{d}\,=\,-\DD^*\cdot A\xi\in L^2(\Omega)$.
For all $\mu\geq \nu  >0$ we have $\phi_\mu=(\mu+\calL)^{-1}\mathfrak{d}$ and
$\phi_\nu=(\nu+\calL)^{-1}\mathfrak{d}$, by the Cauchy-Schwarz inequality,
\begin{eqnarray*}
\expec{|\DD \phi_\mu-\DD\phi_\nu|^2} &\leq& \alpha^{-1} \expec{(\phi_\mu-\phi_\nu)\calL (\phi_\mu-\phi_\nu)} \\
&=& \alpha^{-1} \expec{\phi_\mu\calL\phi_\mu}  -2 \alpha^{-1} \expec{\phi_\mu\calL \phi_\nu} +\alpha^{-1} \expec{\phi_\nu\calL\phi_\nu} \\
&=& \alpha^{-1} \expec{\mathfrak{d}(\mu+\calL)^{-1}\calL(\mu+\calL)^{-1}\mathfrak{d}}  -2 \alpha^{-1} \expec{\mathfrak{d}(\mu+\calL)^{-1}\calL (\nu+\calL)^{-1}\mathfrak{d}} \\
&&+\alpha^{-1} \expec{\mathfrak{d}(\nu+\calL)^{-1}\calL(\nu+\calL)^{-1}\mathfrak{d}} .
\end{eqnarray*} 
By the spectral formula with functions
$$
\Psi(\lambda)\,=\, \frac{\lambda}{(\mu+\lambda)^2},  \frac{\lambda}{(\mu+\lambda)(\nu+\lambda)}, \frac{\lambda}{(\nu+\lambda)^2},
$$
we obtain
\begin{eqnarray}
\expec{|\DD \phi_\mu-\DD\phi_\nu|^2} &\leq& \alpha^{-1}\int_{\R^+} \left(\frac{\lambda}{(\mu+\lambda)^2}-2\frac{\lambda}{(\mu+\lambda)(\nu+\lambda)}+ \frac{\lambda}{(\nu+\lambda)^2}\right) de_{\mathfrak{d}}(\lambda) \nonumber \\
&=& \alpha^{-1}\int_{\R^+} \frac{\lambda(\nu-\mu)^2}{(\mu+\lambda)^2(\nu+\lambda)^2}de_{\mathfrak{d}}(\lambda) \label{eq:ex-sp-calc}\\
&\leq & \alpha^{-1}\int_{\R^+} \frac{\mu^2}{(\mu+\lambda)^2\lambda}de_{\mathfrak{d}}(\lambda)\nonumber ,
\end{eqnarray}
since $0<\nu\leq \mu$. Since the upper bound is independent of $\nu$, we have proved the claim if we can show that it tends to zero as $\mu$ vanishes. This is a consequence of the Lebesgue dominated convergence theorem provided we show that 
\begin{equation}\label{2-1:eq:KV-1}
\int_{\R^+}\frac{1}{\lambda} de_{\mathfrak{d}}(\lambda)\,<\,\infty.
\end{equation}
On the one hand, by the a priori estimate of $\DD \phi_\mu$,
$$
\expec{\phi_\mu \calL \phi_\mu} \,\leq \, \beta \expec{|\DD \phi_\mu|^2} \,\leq \, \beta C.
$$
On the other hand, by the same type of spectral calculus as above, we have
$$
\expec{\phi_\mu \calL \phi_\mu} \,=\,\int_{\R^+}\frac{\lambda}{(\mu+\lambda)^2}de_{\mathfrak{d}}(\lambda) .
$$
Estimate \eqref{2-1:eq:KV-1} then follows from the monotone convergence theorem.
This concludes the proof of Lemma~\ref{2-1:lem:corr-mod-corr}.

\subsection{Random walk in random environment}\label{sec:2-1-1-3}

We now turn our attention to the probabilistic aspects of homogenization. This presentation is informal. It aims at being accessible to non-specialists of probability theory, and at highlighting the inner similarities with the corrector approach of subsection~\ref{sec:2-1-1-2}.

\subsubsection{The continuous-time random walk}
\label{ss:ctrw}

Let the environment $\omega$ be fixed for a while (that is, we have picked a realization of the conductivities $\omega_e \in [\alpha,\beta]$, $e\in \mb B$). The random walk we wish to define, that we will denote by $(X_t)_{t\in \R_+}$, is a random process whose behavior is influenced by the environment. 

To the specialist, it can be defined by saying that it is the Markov process whose transition rates are the $(\omega_e)_{e \in \mb{B}}$. The Markov property means that given any time $t \ge 0$, the behavior of the process after time $t$ depends on its past only through its location at time $t$. In other words, the process ``starts afresh'' at time $t$ given its current location. In order to give a complete description of the process, it thus suffices to describe its behavior over a time interval $[0, t]$, for some $t > 0$, no matter how small. As $t$ tends to $0$, this behavior is given by
\begin{equation}
\label{nonconstr}
\PPo_z \Ll[X_{t} = z' \Rr] = 
\left|
\begin{array}{ll}
t\omega_{z,z'} + o(t) & \text{if } z' \sim z, \\
1-\sum_{y \sim z} t\omega_{z,y} + o(t) & \text{if } z' = z, \\
o(t) & \text{otherwise},
\end{array}
\right.
\end{equation}
where $\PPo_z$ is the probability measure corresponding to the walk started at $z$, that is, $\PPo_z[X_0 = z] = 1$. Equation~\eqref{nonconstr} is what is meant when it is said that $\omega_{z,z'}$ is the jump rate from $z$ to $z'$.

A more constructive way to represent the random walk is as follows. Let the walk be at some site $z \in \Z^d$ at time $t$, and start an ``alarm clock'' that rings after a random time $T$ which follows an exponential distribution of parameter 
\begin{equation}
\label{defpomega}
p_\omega(z) := \sum_{z' \sim z} \omega_{z,z'}.
\end{equation}
This means that for any $s \ge 0$, the probability that $T > s$ is equal to $e^{-p_{\omega(z)} s}$. When the clock rings, the walk chooses to move to one (out of $2d$) neighboring site $z'$ with probability
\begin{equation}
\label{def:leadsto}
p(z\leadsto z') :=
\frac{\omega_{z,z'}}{p_\omega(z)},
\end{equation}
and this choice is made independently of the value of $T$.

Note that by the Markov property, the fact that the walk has not moved during the time interval $[t,t+s]$ should give no information on the time of the next jump.
Only exponential distributions have this memoriless property.

Let us see why thus defined, the random walk satisfies \eqref{nonconstr}. The probability that the clock rings during the time interval $[0,t]$ is
$$
1-e^{-p_\omega(z)t} = p_{\omega(z)} t + o(t).
$$
Since $p_{\omega}(z)$ is bounded by $2d\beta$ uniformly over $z$, the probability that the walk makes two or more jumps is $o(t)$. The probability that it ends up at $z' \sim z$ at time $t$ is thus
$$
\Ll(p_{\omega(z)} t + o(t)\Rr) p(z\leadsto z') - o(t) = \omega_{z,z'} t + o(t),
$$
and the probability that it stays still is indeed as in \eqref{nonconstr}.

\medskip

The link between the random walk and the elliptic operators of the previous subsection is as follows.
We  let $(P_t)_{t \in \R_+}$ be the semi-group associated with the random walk, that is, for any $t \ge 0$ and any bounded function $f:\Z^d \to \R$, we let
\begin{equation*}\label{eq:defi-semi}
P_tf(z)\,=\,\mathbf E_z^\omega[f(X_t)],
\end{equation*}
where $\mathbf E_z^\omega$ denotes the expectation associated to the probability measure $\PPo_z$, under which the random walk starts at $z\in \Z^d$ in the environment $\omega$. This is a semi-group since $X$ has the Markov property. As we now show, the infinitesimal generator of this semi-group is the elliptic operator $-L$ (where $L$ is defined in \eqref{2-1:eq:def-op-space}).
Recall that the infinitesimal generator of a semi-group, applied to $f$, is given by
\begin{equation*}
\lim_{t\to 0} \left(\frac{P_tf-f}{t} \right),
\end{equation*}
for any $f$ for which the limit exists. From the description \eqref{nonconstr}, this limit is easily computed. Indeed,
\begin{eqnarray}
P_tf(z)&=&\mathbf E_z^\omega[f(X_t)]\nonumber \\
&=&(1-p_\omega(z)t)f(z)+t \sum_{z'\sim z} \omega_{(z,z')} f(z')+o(t)\\
&=&f(z)+t \sum_{z'\sim z}\omega_{(z,z')} (f(z')-f(z))+o(t),\label{2-1:eq:bonne-proba2}
\end{eqnarray}
so that by \eqref{2-1:eq:def-op-space2},
\begin{equation}
\label{limsemigrp}
\lim_{t\to 0} \left(\frac{P_tf(z)-f(z)}{t} \right) \,=\, \sum_{z'\sim z}\omega_{(z,z')} (f(z')-f(z))\,=\,-Lf(z),
\end{equation}
and we have identified the infinitesimal generator to be $-L$, as announced.

\medskip

An important feature of this random walk is that the jump rates are symmetric: the probability to go from $z$ to $z'$ in an ``infinitesimal'' amount of time is equal to the probability to go from $z'$ to $z$, as can be seen on \eqref{nonconstr}. This may be rephrased as saying that the counting measure on $\Z^d$ (which puts mass $1$ to every site) is \emph{reversible} for the random walk.

If we were running the random walk on a finite graph instead of $\Z^d$, one could normalize the counting measure to make it into a probability measure, and then interpret reversibility as follows. Start the random walk at a point chosen uniformly at random on the graph, and let it run up to time $t$. Reversibility (of the uniform probability measure) is equivalent to the fact that the distribution of this trajectory is the same as the distribution of the time-reversed path, running from time $t$ to time $0$. In particular, since the starting point is chosen uniformly at random, the distribution of the walk at any time must be the uniform probability measure. Of course, the problem is that on $\Z^d$, it does not make sense to say that we choose the starting point ``uniformly at random'', or in other words, we cannot normalize the counting measure into a probabability measure.

In the previous section, we moved from the operator $L$ to its ``environmental'' version, the operator $\L$. This takes a very enlightening probabilistic meaning. Instead of considering the random walk itself, we may consider the \emph{environment viewed by the particle}, which is the random process defined as
$$
t \mapsto \omega(t) := \theta_{X(t)}\omega,
$$
where $(\theta_x)_{x \in \Z^d}$ are the translations introduced in Definition~\ref{2-1:defi:stat}. One can convince oneself that $(\omega(t))_{t \in \R_+}$ is a Markov process, and the important point is that its infinitesimal generator is precisely $-\L$. A little computation (see for instance \cite[Proposition~3.1]{these}) shows that the reversibility of the counting measure for the random walk translates into the reversibility of the measure $\P$ for the process $(\omega(t))_{t \in \R_+}$. This reversibility implies that $(\omega(t))_{t \in \R_+}$ is stationary if we take the initial environment according to $\P$ and then let the process run, that is, under the measure
\begin{equation}\label{2-1:eq:def-pdt-measure}
\P_0 :=\P \mathbf P_0^\omega,
\end{equation}
the so-called \emph{annealed} measure.
Stationarity means that for any positive integer $k$ and any $s_1,\ldots, s_k \ge 0$, the distribution of the vector $(\omega(s_1+t),\ldots, \omega(s_k+t))$ does not depend on $t \ge 0$. In particular, under the measure $\P_0$, the distribution of $\omega(t)$ is the measure ${\P}$ for any $t \ge 0$. 

To conclude, we argue that from the reversibility follows the fact that the operator $\L$ is self-adjoint in $L^2(\Omega)$ --- which we have already seen more directlty in the previous subsection. Indeed, the invariance under time reversal implies that, under $\P_0$ and for any $t \ge 0$, the vectors $(\omega(0), \omega(t))$ and $(\omega(t), \omega(0))$ have the same distribution. For any bounded functions $f,g :\Omega \to \R$, we thus have
\begin{equation}
\label{presqselfadj}
\E_0\Ll[ f(\om(0))\ g(\om(t)) \Rr] = \E_0\Ll[ f(\om(t))\ g(\om(0)) \Rr],
\end{equation}
where we write $\E_0$ for the expectation associated to $\P_0$. Letting $(\mcl{P}_t)_{t \in \R_+}$ be the semi-group associated to $-\L$, we arrive at
\begin{eqnarray*}
\E_0\Ll[ f(\om(0))\ g(\om(t)) \Rr] & = & \expec{ \EEo_0 \Ll[ f(\om(0)) \ g(\om(t)) \Rr] } \\
& = & \expec{  f(\om) \ \EEo_0 \Ll[ g(\om(t)) \Rr] } = \expec{ f \ \mcl{P}_t g},
\end{eqnarray*}
and \eqref{presqselfadj} thus reads
$$
\expec{ f \ \mcl{P}_t g} = \expec{ \mcl{P}_t f \ g}.
$$
This shows that $\mcl{P}_t$ is self-adjoint for any $t \ge 0$. Passing to the limit as in \eqref{limsemigrp}, we obtain the self-adjointness of $\L$ itself. In fact, reversibility of a Markov process and self-adjointness of its infinitesimal generator are two sides of the same coin.

\subsubsection{Central limit theorem for the random walk}
\label{ss:clt}

The aim of this paragraph is to sketch a probabilistic argument justifying the following result.
\begin{thm}[\cite{Kipnis-Varadhan-86}]
\label{t:kv}
Under the measure $\P_0$ and as $\eps$ tends to $0$, the rescaled random walk $X^{(\eps)} := (\sqrt{\e}X_{t/\e})_{t \in \R_+}$ converges in distribution (for the Skorokhod topology) to a Brownian motion with covariance matrix $2A_\ho$, where $A_\ho$ is as in \eqref{2-1:eq:hom-coeff}. In other words, for any bounded continuous functional $F$ on the space of cadlag functions, one has
\begin{equation}
\label{convergegen}
\E_0\Ll[ F(X^{(\eps)}) \Rr] \xrightarrow[\eps \to 0]{} E[F(B)],
\end{equation}
where $B$ is a Brownian motion started at the origin and with covariance matrix $2 \Ah$, and $E$ denotes averaging over $B$.  Moreover, for any $\xi \in \R^d$, one has
\begin{equation}
\label{convergesquare}
 t^{-1} \E_0\Ll[ \Ll(\xi \cdot X_t\Rr)^2 \Rr] \xrightarrow[t \to + \infty]{} 2 \xi \cdot \Ah \xi.
\end{equation}
\end{thm}
Note that the convergence of the rescaled square displacement in \eqref{convergesquare} does not follow from \eqref{convergegen}, since the square function is not bounded.

From now on, we focus on the one-dimensional projections of $X$. As in subsection~\ref{sec:2-1-1-2}, we let $\xi$ be a fixed vector of $\R^d$. The idea is to decompose $\xi \cdot X_t$ as 
\begin{equation}
\label{decomp}
\xi \cdot X_t = M_t + R_t,
\end{equation}
where $(M_t)_{t \ge 0}$ is a martingale and $R_t$ is a (hopefully small) remainder. 

We recall that $(M_t)_{t \ge 0}$ is a martingale under $\mathbf E_0^\omega$ if for all $t\geq 0$ and $s\geq 0$,
\begin{equation}\label{2-1:eq:am-I-martingale?}
\mathbf E^\omega_0[M_{t+s}|\calF_t]\,=\,M_t,
\end{equation}
where $\calF_t$ is the $\sigma$-algebra generated by $\{M_\tau,\tau\in [0,t]\}$. (If $M$ was a gambler's money, one should say that he is playing a fair game when \eqref{2-1:eq:am-I-martingale?} holds, since knowing the history up to time $t$, his expected amount of money at time $t+s$ is the amount he has at time $t$). 

We look for a martingale of the form $M_t=\chi^\omega(X_t)$ for some function $\chi^\omega$. By the Markov property of $X$, $(M_t)_{t \ge 0}$ of this form is a martingale if and only if, for any $z \in \Z^d$ and any $t \ge 0$,
\begin{equation}\label{eq:martingale}
\mathbf E^\omega_z[\chi^\omega(X_t)]\,=\,\chi^\omega(z),
\end{equation}
From \eqref{nonconstr}, we learn that
$$
\mathbf E^\omega_z[\chi^\omega(X_t)] = \chi^\om(z) - s L\chi^\omega(z) + O(s^2).
$$
Hence, a necessary condition is that $L\chi^\om(z) = 0$, and in fact, this condition is also sufficient. Keeping in mind that we also want the remainder term to be small, we would like $\chi^\omega$ to be a perturbation of the $z \mapsto \xi \cdot z$, so that a right choice for $\chi^\omega$ should be
\begin{equation*}
\chi^\omega(z)\,=\,\xi \cdot z+\phi(z,\omega)
,
\end{equation*}
where $\phi$ is the corrector of Lemma~\ref{2-1:lem:corr} (compare equation~\eqref{2-1:eq:corr-sto} to $L\chi^\om(z) = 0$). The link between the corrector equation and the RWRE appears precisely there. We thus define
$$
M_t = \xi \cdot X_t + \phi(X_t,\omega) \quad \text{and} \quad R_t = -\phi(X_t,\omega).
$$
We now argue that the martingale $M$ has stationary increments under the measure $\P_0$. For simplicity, let us just see that the distribution of $M_{t+s} - M_t$ under $\P_0$ does not depend on~$s$ (instead of considering vectors of increments). The argument relies on the stationarity of the process of the environment viewed by the particle seen above. Given this stationarity, it suffices to see that $M_{t+s} - M_t$ can be written as a function of $(\omega(t),\omega(t+s))$ only. This is possible because the increment $M_{t+s} - M_t$ depends only on~$\nabla \phi$, which is a function of $\omega$ only (contrary to $\phi$ itself which is a priori a function of $x$ and $\omega$).

\medskip

Martingales are interesting for our purpose since they are ``Brownian motions in disguise''. To make this idea more precise, let us point out that any one-dimensional continuous martingale can be represented as a time-change of Brownian motion (this is the Dubins-Schwarz theorem, see \cite[Theorem~V.1.6]{ry}). If $(B_t)_{t \ge 0}$ is a Brownian motion, a time-change of it is for instance $M_t = B_{t^{7}}$. Note that in this example, the time-change can be recovered by computing $E[M_t^2] = t^7$.
Here, the martingale we consider has jumps (since $X$ itself has), which complicates the matter a little, but let us keep this under the rug. Intuitively, in order to justify the convergence to a Brownian motion, we should show that the underlying time-change grows linearly at infinity. One can check that two increments of a martingale over disjoint time intervals are always orthogonal in $L^2$ (provided integration is possible). In our case, since $M_t$ has stationary increments, it thus follows that letting $\sigma^2(\xi) = \E_0[M_1^2]$, we have
\begin{equation}
\label{meanlinear}
\E_0[M_t^2] = \sigma^2(\xi) \ t,
\end{equation}
so we are in good shape (i.e.\ on a heuristic level, it indicates that the underlying time-change indeed grows linearly). Letting $f:z\mapsto (x\cdot \xi+\phi(x,\omega)\cdot \xi)^2$ and using \eqref{nonconstr}, one can write
\begin{eqnarray*}
\lefteqn{\mathbf E^\omega_0[(M_t\cdot \xi)^2]\,=\,\mathbf E^\omega_0[f(X_t)]}\\
&= & f(0)+t \sum_{z'\sim 0}\omega_{(z,z')}\big(z'\cdot \xi+\phi(z',\omega)\cdot \xi \big)^2+o(t)\\
&=&\underbrace{2t (\xi+\nabla \phi(0,\omega)\cdot \xi)\cdot A(0,\omega) (\xi+\nabla \phi(0,\omega)\cdot \xi)}_{=: t\  v(\omega)} + o(t),
\end{eqnarray*}
since $\phi(0,\omega)=0$. From the definition of $\Ah$ in \eqref{2-1:eq:hom-coeff}, we thus get 
\begin{equation}
\label{identifcv}
\sigma^2(\xi) = \expec{v} = 2 \xi \cdot \Ah \xi.
\end{equation}
Provided we can show that the remainder is negligible, this already justifies \eqref{convergesquare}. 
In order to prove that $(\sqrt{\eps} M_{\eps^{-1} t})_{t \ge 0}$ converges to a Brownian motion of variance $\sigma^2(\xi)$ as $\eps$ tends to~$0$, knowing \eqref{meanlinear} is however not sufficient: one does not recover all the information about the time-change by computing $\E_0[M_t^2]$ alone. This can be understood from the fact that in the Dubins-Schwarz theorem, the time-change that appears can be random itself, and is in fact the quadratic variation of the martingale. We will not go into explaining what the quadratic variation of a martingale is in general, but simply state that in our case, it is given by
$$
V_t = \int_0^t v(\omega(s)) \ \d s,
$$
and what we should prove is thus that
\begin{equation}
\label{ergodicthm}
t^{-1} V_t = t^{-1} \int_0^t v(\omega(s)) \ \d s\xrightarrow[t \to + \infty]{\text{a.s.}} \sigma^2(\xi) = 2 \xi \cdot \Ah \xi.
\end{equation}
One can show that the process $(\omega(t))_{t \ge 0}$ is ergodic for the measure $\P_0$ (see for instance \cite[Proposition~3.1]{these}), and thus the convergence in \eqref{ergodicthm} is a consequence of the ergodic theorem. Note in passing this surprising fact that the proof of a central limit theorem was finally reduced to a law of large numbers type of statement.

\medskip

In order to obtain a central limit theorem for $\xi \cdot X$ itself (instead of the martingale part), we need to argue that the remainder term is small. In order to do so, and following \cite{Kipnis-Varadhan-86}, we will rely on spectral theory.
We recall that the operator $\calL$ of Definition \ref{2-1:defi:operator-sto} being a bounded non-negative self-adjoint operator on $L^2(\Omega)$, it admits a spectral decomposition in $L^2(\Omega)$.
Moreover, for any $g\in L^2(\Omega)$ the characterizing property of the spectral measure~$e_g$ previously defined is that for any continuous function $\Psi:[0,+\infty)\to \R_+$, one has
\begin{equation*}
\expec{g \ \Psi(\calL)g}\,=\,\int_{\R^+}\Psi(\lambda) \ \d e_g(\lambda). 
\end{equation*}
We now argue that
\begin{equation}\label{2-1:eq:KV-2}
\frac{1}{t}\overline{\mb E}[R_t^2]\,=\,2\int_{\R^+}\frac{1-e^{-t\lambda}}{t\lambda^2} de_{\mathfrak{d}}(\lambda) \xrightarrow[t \to + \infty]{} 0,
\end{equation}
where $\mathfrak{d}$ is the local drift in the direction $\xi$, that is, 
\begin{equation}
\label{defmfkd}
\mfk{d} = - \nabla^* \cdot A(0,\omega) \xi = -\DD^* \cdot A \xi = \sum_{z \sim 0} \omega_{0,z} \xi \cdot z.
\end{equation}
We start by showing the equality in \eqref{2-1:eq:KV-2}, and will later show that the spectral integral tends to $0$ as $t$ tends to infinity. Recall that
$$
R_t = -\phi(X_t,\om) = \phi(0,\om)-\phi(X_t,\om),
$$
and that $\phi(0,\om)-\phi(x,\om)$ can be obtained as the limit of $\phi_\mu(0,\om)-\phi_\mu(x,\om)$. To make them easier, we do the computations below as if $\phi$ was a stationary $\phi_\mu$. The argument can be made rigorous through spectral analysis as in \cite[Theorem~8.1]{Mourrat-10}, or using Lemma~\ref{2-1:lem:corr-mod-corr}. We expand the square:
$$
\E_0[R_t^2] = \E_0[(\phi(0,\om))^2] + \E_0[(\phi(X_t,\om))^2] - 2 \E_0[\phi(0,\om)\phi(X_t,\om)].
$$
Note that, by the simplifying assumption, 
$$
\E_0[(\phi(X_t,\om))^2] = \E_0[(\phi(0,\theta_{X_t}\ \om))^2] = \E_0[(\phi(0,\om(t)))^2].
$$
Since $(\om(t))_{t \in \R_+}$ is stationary under $\P_0$, the last expectation is equal to $\E_0[(\phi(0,\om))^2] = \expec{\phi^2}$. Since $\L \phi = \mfk{d}$, we have
$$
\expec{ \phi^2} =\expec{\L^{-1} \mfk{d} \ \L^{-1} \mfk{d}} =  \expec{\mfk{d} \ \L^{-2} \mfk{d}} = \int \lambda^{-2} \ \d e_\mfk{d}(\lambda).
$$
For the cross-product, 
\begin{equation}
\label{I-semig}
\E_0[\phi(0,\om)\phi(X_t,\om)] = \expec{ \phi(0,\omega) \ \EEo_0\Ll[ \phi(0,\om(t)) \Rr]  },
\end{equation}
and $\EEo_0\Ll[ \phi(0,\om(t)) \Rr]$ is the image $\phi$ by the semi-group associated to $-\L$, that is, $e^{-t\L}$. Using also  the fact that $\L \phi = \mfk{d}$, we can rewrite the r.h.s.\ of \eqref{I-semig} as
$$
\expec{\phi \ e^{-t\L} \phi} = \expec{\L^{-1} \mfk{d} \ e^{-t\L} \L^{-1} \mfk{d}} = \int \lambda^{-2} e^{-t\lambda} \ \d e_\mfk{d}(\lambda),
$$
and this justifies the equality \eqref{2-1:eq:KV-2}. The fact that the spectral integral in \eqref{2-1:eq:KV-2} tends to zero follows from the dominated convergence theorem and \eqref{2-1:eq:KV-1}.
\medskip

Before concluding this section, we point out several differences between the argument as we presented it and the original one from \cite{Kipnis-Varadhan-86}. Here, we used the existence of the corrector, borrowed from subsection~\ref{sec:2-1-1-2}, to construct the martingale. In \cite{Kipnis-Varadhan-86}, the martingale is constructed directly, by considering the martingale
$$
\xi \cdot X_t + \phi_\mu(\omega(t)) - \phi_\mu(\omega(0)) - \mu \int_0^t \phi_\mu(\omega(s)) \ \d s,
$$
and showing that for each fixed $t$, it is a Cauchy sequence in $L^2(\P_0)$ (and thus converges) as $\mu$ tends to $0$. This is achieved through spectral analysis, using the estimate \eqref{2-1:eq:KV-1}. This estimate is obtained through a general argument of (anti-) symmetry (see the proof of \cite[Theorem~4.1]{Kipnis-Varadhan-86}), which has been systematized by \cite{DeMasi-Ferrari-89}. Another difference is that the arguments of \cite{Kipnis-Varadhan-86} are developped for general reversible Markov processes.

\section{Numerical approximation of the homogenized coefficients using the corrector equation}
\label{s:num-corr}
\setcounter{equation}{0}

\subsection{General approach}

Let $\xi\in \R^d$ be a fixed unit vector.
In order to approximate $A_\ho$ using the corrector $\phi$,  we first replace the expectation in \eqref{2-1:eq:hom-coeff} by a spatial average appealing to ergodicity: 
almost surely
\begin{equation}\label{eq:as-fo-app}
\xi\cdot A_\ho \xi\,=\, \lim_{N \to \infty} \fint_{Q_N}(\xi+\nabla \phi)\cdot A(\xi+\nabla \phi),
\end{equation}
where $Q_N=[0,N)^d$ and $\fint_{Q_N}:=N^{-d}\sum_{Q_N\cap \Z^d}$.

Recall that $\phi$ is the solution to \eqref{2-1:eq:corr-sto}, which is a problem posed on $\Z^d$.
In order to turn \eqref{eq:as-fo-app} into a practical formula, one needs a computable approximation $\tilde \phi_N$ of $\phi$ on $Q_N$,
which would allow us to define 
an approximation $A_N$ of $A_\ho$ by
\begin{equation*}
\xi\cdot A_N \xi\,=\, \fint_{Q_N}(\xi+\nabla \tilde \phi_N)\cdot A(\xi+\nabla \tilde \phi_N).
\end{equation*}
Although $A_\ho$ is deterministic, this approximation $A_N$ is a random variable.
Provided the chosen approximation $\tilde \phi_N$ is ``consistent'', 
\begin{equation*}
\lim_{N\to \infty}\xi\cdot A_N \xi\,=\,\xi\cdot A_\ho\xi
\end{equation*}
almost surely.

\medskip

The starting point for a quantitative convergence analysis is the following identity:
\begin{equation}\label{eq:starting-point}
\expec{(\xi\cdot A_N\xi-\xi\cdot A_\ho\xi)^2} \,=\,\var{\xi\cdot A_N\xi}+(\xi\cdot \expec{A_N}\xi-\xi\cdot A_\ho\xi)^2.
\end{equation}
The square root of the first term of the r.~h.~s. is called the \emph{random error} or \emph{statistical error}. It measures the fluctuations of $A_N$ around its expectation.
The square root of second term of the r.~h.~s. is the \emph{systematic error}. It measures the fact that $\nabla \tilde \phi_N$ is only an approximation of $\nabla \phi$ on $Q_N$.
In order to give quantitative estimates of these errors,
we need to make assumptions on the statistics of $A$, and shall assume in the rest of this section that the entries of $A$ have finite correlation length (that is, there exists $\mathscr{L}_c$
such that $a(x,x+\ee_i)$ and $a(z,z+\ee_j)$ are independent if $|x-z|\geq \mathscr{L}_c$, for all $i,j\in \{1,\dots,d\}$), which sligthly generalises the i.i.d.\ case and is enough for our purposes.

\medskip

In the following subsection, we introduce three approximation methods based on the above strategy, and recall the known (or expected) convergence rates
for both errors.
The third subsection is dedicated to a numerical study which completes the analysis in two respects: it allows to confirm/infirm numerically the expected convergence rates (if they are not precisely known), and it allows to make explicit the prefactors.

\subsection{Methods and theoretical analysis}

\subsubsection{Homogeneous Dirichlet boundary conditions}

The simplest approximation of the corrector $\phi$ on $Q_N$ for $N\geq 1$ is given by the unique solution $\phi_N\in L^2(Q_N\cap \Z^d)$ to 
\begin{equation*}
-\nabla \cdot A(\xi+\nabla \phi_N)\,=\,0 \qquad \text{ in }Q_N,
\end{equation*}
completed by the boundary conditions $\phi_N(x)=0$ for all $x\in \Z^d\setminus Q_N$.

\medskip

The approximation $A_N$ of $A_\ho$ is then defined by
\begin{equation*}
\xi\cdot A_N\xi\,:=\,\fint_{Q_N} (\xi+\nabla \phi_N)\cdot A(\xi+\nabla \phi_N).
\end{equation*}
As a direct consequence of homogenization, we have the almost sure convergence:
\begin{equation*}
\lim_{N\to \infty} |A_N-A_\ho|\,=\,0.
\end{equation*}
In terms of convergence rates, the starting point is identity \eqref{eq:starting-point}.
Although we do not have a complete proof, we expect the random error to scale as the central limit theorem, that is:
\begin{equation}\label{eq:dir-random}
\var{\xi \cdot A_N\xi}^{1/2}\,\sim \, N^{-d/2},
\end{equation}
and the systematic error to scale as a surface effect (the corrector is perturbed on the boundary):
\begin{equation}\label{eq:dir-syst}
|\xi \cdot \expec{A_N}\xi-\xi\cdot A_\ho\xi|\,\sim \, N^{-1}.
\end{equation}
Let us give an informal argument for \eqref{eq:dir-random} in the case of ellipticity ratio $\beta/\alpha$ close to $1$, that is, for $A$ a perturbation of $\Id$. Then, at first order, the equation for $\phi_N$ takes the form
$$
-\triangle \phi_N \,=\,\nabla \cdot A \xi,
$$
completed by homogeneous Dirichlet boundary conditions on $\partial Q_N$. Let $G$ denote the Green's function of the discrete Laplace equation on $Q_N$ with homogeneous Dirichlet boundary conditions.
Then, by the Green representation formula,
\begin{eqnarray}\label{eq:motiv-var}
\phi_N(x)&=&\int_{Q_N} \nabla_y G(x,y) \cdot A(y) \xi dy,\\
\nabla \phi_N(x)&=&\int_{Q_N} \nabla_x \nabla_y G(x,y) \cdot A(y) \xi dy.\label{eq:motiv-var2}
\end{eqnarray}
We then appeal to the following spectral gap inequality (which is at the core \cite{Gloria-Otto-09}): for any function $X$ of a finite number of the i.i.d.\ random variables $\omega_e$, we have:
\begin{equation}\label{eq:sge}
\var{X} \,\leq \, \sum_e \expec{\sup_{\omega_e} \left|\frac{\partial X}{\partial \omega_e}\right|^2}\var{\omega},
\end{equation}
where the supremum is taken w.~r.~t. the variable $\omega_e$, and $\var{\omega}$ is the variance of the i.i.d.\ conductances.
We shall apply this inequality to $X=\xi\cdot A_N\xi$.
We first note that the weak formulation of the equation yields
$$
X=\xi\cdot A_N\xi\,=\,\fint_{Q_N} \xi\cdot A(\xi+\nabla \phi_N),
$$
so that by \eqref{eq:motiv-var2}, we have for all $e=(z,z+e_i)$
$$
\frac{\partial X}{\partial \omega_e}\,=\,\frac{1}{N^d}\xi_i(\xi_i+\nabla_i \phi_N(z))+\frac{1}{N^d} \int_{Q_N} \xi \cdot A(x)\nabla_x \nabla_{z_i} G(x,z)  \xi_i dx,
$$
which, by symmetry of $A$ and of $G$, we rewrite as 
$$
\frac{\partial X}{\partial \omega_e}\,=\,\frac{1}{N^d}\xi_i(\xi_i+\nabla_i \phi_N(z))+\frac{1}{N^d} \xi_i \nabla_{z_i} \int_{Q_N} \nabla_x G(z,x) \cdot A(x) \xi  dx.
$$
Using \eqref{eq:motiv-var}, this turns into 
$$
\frac{\partial X}{\partial \omega_e}\,=\,\frac{1}{N^d}\xi_i(\xi_i+2\nabla_i \phi_N(z)).
$$
It remains to take the supremum of this quantity w.~r.~t.  $\omega_e$. In view of \eqref{eq:motiv-var}, we have
for all $j$
$$
\osc{\omega_e}{\nabla_j \phi(z)}=\sup_{\omega_e} \nabla_j \phi(z)
- \inf _{\omega_e} \nabla_j \phi(z)\,\leq\, |\nabla_{x_j} \nabla_{y_i} G(x,z)|,
$$
which is bounded by a universal constant (in the discrete case the Green function is bounded).
We thus have the desired variance estimate:
\begin{eqnarray*}
\var{X} &\leq & \frac{1}{N^{2d}}\sum_e \expec{C(1+|\nabla_i \phi_N(z)|^2)} \var{\omega} \\
&=&\frac{1}{N^{2d}}\expec{\sum_e C(1+|\nabla_i \phi_N(z)|^2)} \var{\omega} \\
&=& \frac{1}{N^d}   \, \frac{1}{N^d}  \expec{ d C N^d+ C\|\nabla \phi_N\|_{L^2(Q_N)}^2}  \var{\omega} \\
&\lesssim & \frac{1}{N^d} ,
\end{eqnarray*}
since an elementary deterministic estimate yields
$$
\int_{Q_N}|\nabla \phi_N|^2 \,\leq \, \beta^2 |Q_N| = \beta^2 N^d.
$$
The difficulty in the case of general ellipticity ratio is to treat the dependence of the Green's function with respect to $\omega_e$, which yields additional terms of the form $|\nabla \phi_N|^4$ which we do not control a priori (see \cite{Gloria-Neukamm-Otto-a} in the case of periodic approximations), except for ellipticity ratio close to 1 (so that Meyers' estimate yields sufficient additional integrability of $\nabla \phi_N$).

\medskip

Note that the convergence rate of the random error \eqref{eq:dir-random} is expected to depend on the dimension whereas the convergence rate of the systematic
error \eqref{eq:dir-syst} does not. 
The combination of these (conjectured) estimates would yield the following convergence rate in any dimension
$$
\expec{(\xi\cdot A_N\xi-\xi\cdot A_\ho\xi)^2}\,\sim \,N^{-2}.
$$

\subsubsection{Regularized corrector and filtering}

In \cite{Gloria-Otto-09,Gloria-Otto-09b,Gloria-10}, the following strategy was used.
Instead of considering an approximation of the corrector $\phi$, we consider an approximation of the
regularized corrector $\phi_\mu$ of Lemma~\ref{2-1:lem:mod-corr}.
The advantage of the regularized corrector is that the Green's function associated with the operator $\mu-\nabla^*\cdot A\nabla$
decays exponentially in terms of the distance measured in units of $\mu^{-1/2}$.
Let $N>0$. We denote by $\phi_{\mu,N}$ the unique solution in $L^2(Q_N\cap \Z^d)$ to 
\begin{equation*}
\mu \phi_{\mu,N}-\nabla \cdot A(\xi+\nabla \phi_{\mu,N})\,=\,0 \qquad \text{ in }Q_N,
\end{equation*}
completed by the boundary conditions $\phi_{\mu,N}(x)=0$ for all $x\in \Z^d\setminus Q_N$.

To define an approximation of $A_\ho$, we introduce for all $L\in \N$ an averaging mask $\eta_L:\Z^d\to \R^+$, with support in $Q_L$, and such that
$\int_{Q_L} \eta_L=1$ and $\sup_{\Z^d} |\nabla \eta_L| \lesssim L^{d+1}$.
For all $N\geq L>0$ and $\mu>0$, we then set:
\begin{equation}\label{def:Azero}
\xi\cdot A_{\mu,N,L}\xi \,:=\,\int_{Q_L} (\xi+\nabla  \phi_{\mu,N})\cdot A  (\xi+\nabla  \phi_{\mu,N})\eta_L.
\end{equation}
We have the following almost sure convergence:
$$
\lim_{\mu\to 0} \lim_{N\geq L\to \infty}|A_{\mu,N,L}-A_\ho|\,=\,0.
$$

\medskip

Following the general approach, we have using the triangle inequality, and expanding the square:
\begin{equation}
\expec{(\xi\cdot(A_{\mu,N,L}-A_{\ho})\xi)^2}^{1/2} \leq \, \expec{(\xi\cdot(A_{\mu,N,L}-A_{\mu,L})\xi)^2}^{1/2}+\var{A_{\mu,L}}^{1/2}+|A_\mu-A_\ho|,\label{eq:triangle-zero}
\end{equation}
where:
\begin{eqnarray}
\xi\cdot A_{\mu,L}\xi&:=&\int_{Q_L} (\xi+\nabla  \phi_{\mu})\cdot A  (\xi+\nabla  \phi_{\mu})\eta_L,\nonumber\\
\xi\cdot A_{\mu}\xi&:=&\expec{(\xi+\nabla  \phi_{\mu})\cdot A  (\xi+\nabla  \phi_{\mu})}.\label{eq:Amu}
\end{eqnarray}
The first term is the error due to the fact that we replace $\phi_\mu$ by the computable approximation $\phi_{\mu,N}$ on $Q_L$.
As shown in \cite{Gloria-10} the function $\phi_{\mu,N}$ is a good approximation of $\phi_\mu$ in the sense
that we have the following deterministic estimate: there exists $c>0$ such that for all $0<L\leq N$, we have almost surely
\begin{equation}\label{eq:syst-det}
|A_{\mu,N,L}-A_{\mu,L}| \,\lesssim \, \mu^{-3/4}\left(\frac{N}{L}\right)^{d/2}\left(\frac{N}{N-L}\right)^{d-1/2} \exp\left(-c \sqrt{\mu} (N-L) \right) .
\end{equation}
This estimate essentially follows from the exponential decay of the Green's function measured in units of $\mu^{-1/2}$.
Hence, provided $\sqrt{\mu} (N-L) \gg 1$, the error between $A_{\mu,N,L}$ and $A_{\mu,L}$ is negligible.

Within the assumption of finite correlation length, it is shown in \cite{Gloria-Otto-09} that there exists $q>0$ (depending only on $\alpha$ and $\beta$) such that 
\begin{equation}\label{eq:0dir-random}
\var{\xi \cdot A_{\mu,L}\xi}\,\lesssim \, \left|
\begin{array}{rcl}
d=2&:&L^{-2}\ln^q \mu+\mu^2 \ln^q\mu,\\
d>2&:&L^{-d}+\mu^2L^{2-d},
\end{array}
\right.
\end{equation}
that is the central limit theorem scaling (up to a logarithmic correction in dimension 2) provided $\mu L \lesssim 1$. (We expect this convergence rate to be optimal in general.)

In \cite{Gloria-Otto-09b} (see also \cite{Gloria-Neukamm-Otto-a} for $d=2$), it is proved that  
\begin{equation}\label{eq:0dir-syst}
|\xi\cdot A_{\mu}\xi-\xi\cdot A_\ho\xi|\,\lesssim \, 
\left|
\begin{array}{rcl}
d=2&:&\mu,\\
d=3&:&\mu^{3/2},\\
d=4&:&\mu^{2}|\ln\mu|,\\
d>4&:&\mu^{2}.
\end{array}
\right.
\end{equation}
Note that this scaling depends on the dimension and saturates at dimension $d=4$.
The proof of this estimate is interesting because similar arguments are used to analyze methods based on RWRE. It relies on an observation by one of us in \cite{Mourrat-10}, and spectral calculus.
Recall that for fixed unit vector $\xi \in \R^d$ we denote by $\mathfrak{d}=\DD^* \cdot A \xi$ the local drift and by $e_{\mathfrak{d}}$ the projection of the spectral measure onto $\mathfrak{d}$. By definition,
\begin{eqnarray*}
\xi\cdot A_\mu\xi-\xi\cdot A_\ho\xi &=& \expec{(\xi+\DD \phi_\mu)\cdot A(\xi+\DD\phi_\mu)-(\xi+\DD \phi)\cdot A(\xi+\DD\phi)}\\
&=&\expec{(\DD \phi_\mu-\DD \phi)\cdot A(\xi+\DD\phi)+(\xi+\DD \phi_\mu)\cdot A(\DD \phi_\mu-\DD \phi)}.
\end{eqnarray*}
Using the weak form of the corrector equation \eqref{eq:corr-eq-ant1}, one has for all $\psi \in L^2(\Omega)$,
$$
\expec{\DD \psi\cdot A(\xi+\DD\phi)} \,=\,0.
$$
Taking $\psi=\phi_\mu-\phi_\nu$ for $\nu>0$, and using the symmetry of $A$, one obtains
$$
\expec{(\DD \phi_\mu-\DD \phi_\nu)\cdot A(\xi+\DD\phi)}\,=\,0\,=\,-\expec{(\xi+\DD\phi)\cdot A (\DD \phi_\mu-\DD \phi_\nu)}.
$$
Taking the limit $\nu\to 0$ and using the strong convergence of $\DD \phi_\nu$ to $\DD \phi$ given by Lemma~\ref{2-1:lem:corr-mod-corr}, this finally yields
$$
\xi\cdot A_\mu\xi-\xi\cdot A_\ho\xi \,=\,\expec{(\DD \phi_\mu-\DD \phi)\cdot A(\DD \phi_\mu-\DD \phi)}.
$$
Hence, by spectral calculus and taking the limit $\nu \to 0$ in \eqref{eq:ex-sp-calc}, we obtain
$$
0\,\leq \, \xi\cdot A_\mu\xi-\xi\cdot A_\ho\xi  \,=\,  \mu^2 \int_{\R^+} \frac{1}{(\mu+\lambda)^2\lambda}de_{\mathfrak{d}}(\lambda).
$$
This immediately implies that $\xi\cdot A_\mu\xi-\xi\cdot A_\ho\xi \sim \mu^2$ if $\lambda \mapsto \lambda^{-3}$ is $de_{\mathfrak{d}}$-integrable on $\R^+$.
More precisely, what matters in this spectral integral is the behavior at the bottom of the spectrum. We cannot expect a spectral gap (which would hold in the periodic case since there is a Poincar\'e's inequality on the torus) but we expect the bottom of the spectrum to be sufficiently thin to yield the result.
It is convenient to rewrite the spectral integral as
$$
\int_{\R^+} \frac{1}{(\mu+\lambda)^2\lambda}de_{\mathfrak{d}}(\lambda) \,\leq \, \int_{0}^1 \frac{1}{(\mu+\lambda)^2\lambda}de_{\mathfrak{d}}(\lambda) +\int_{\R^+} \frac{1}{\lambda}de_{\mathfrak{d}}(\lambda) .
$$
By \eqref{2-1:eq:KV-1} the second term of the r.\ h.\ s.\  is bounded, and an elementary calculation shows that \eqref{eq:0dir-syst} is a consequence of the following optimal estimate of the so-called ``spectral exponents": for all $0<\nu\leq 1$ and $d\geq 2$,
\begin{equation}\label{eq:optima-sp-exp}
\int_0^\nu de_{\mathfrak{d}}(\lambda)\,\lesssim \, \nu^{d/2+1}.
\end{equation}
Suboptimal bounds with exponents $d/2-2$ were first obtained for $d$ large in \cite[Theorem~2.4]{Mourrat-10} using probabilistic arguments, optimal bounds up to dimension $4$ (and up to some logarithm in dimension 2) in \cite{Gloria-Otto-09b} using the spectral gap estimate \eqref{eq:sge} and elliptic regularity theory, up to dimension $6$ in \cite{Gloria-Mourrat-10a} by pushing forward the method of \cite{Gloria-Otto-09b} (see \cite{Gloria-Mourrat-10b} for a generalization of this strategy which would yield the optimal exponents
for all $d>2$), and in any dimension in \cite{Gloria-Neukamm-Otto-a} using the spectral gap estimate and parabolic regularity theory to obtain first bounds on the semi-group which turn after integration in time into bounds on the spectral exponents.

\medskip

Altogether, if we take  $N-L \sim N\sim L$, and $\mu\sim N^{-\gamma}$ with $1 \leq \gamma<2$ (so that \eqref{eq:syst-det} is of infinite order in $N$), the combination of the three estimates \eqref{eq:syst-det}, \eqref{eq:0dir-random}, and \eqref{eq:0dir-syst}  yields:
\begin{equation}\label{eq:error-zero}
\expec{|A_{\mu,N,L}-A_{\ho}|^2} \,\lesssim \,  
\left|
\begin{array}{rcl}
d=2&:&N^{-2}\ln^q N,\\
d=3&:&N^{-3},\\
d=4&:&N^{-4}+N^{-4\gamma}\ln^2 N,\\
d>4&:&N^{-4},
\end{array}
\right.
\end{equation}
which yields the central limit theorem scaling up to dimension 4 at least.

\medskip

In dimension $d>2$, since the fluctuations of $A_{\mu,L}$ have the scaling of the central limit theorem \eqref{eq:0dir-random}, one may wonder whether the distribution of
$L^{d/2}(\xi\cdot A_{\mu,L}\xi-\xi\cdot \expec{A_{\mu,L}}\xi)$ converges in law to a Gaussian random variable.
Related results were obtained by Nolen \cite{Nolen-11}, Rossignol \cite{Rossignol-12}, and Biskup, Salvi, and Wolff \cite{Biskup-Salvi-Wolff-12}.

\subsubsection{Periodization method}
\label{ss:period}
The periodization method is a widely used method to approximate homogenized coefficients, which consists in periodizing the random medium,
see \cite{Owhadi-03,E-Yue-07}.
What is less known is that there may be several ways to periodize a random medium:
\begin{itemize}
\item the periodization in space,
\item the periodization in law.
\end{itemize}
Both periodization methods coincide in the specific case of i.i.d.\ conductances. The periodization in law is implicitly used in the very nice contribution \cite{KFGMJ-03}.

\medskip

Let $A$ be a random conductivity function with finite correlation length.
Let us begin with the periodization in space.
It consists in approximating the corrector $\phi$ on $Q_N$ by the $Q_N$-periodic solution $\phi_N^\spa$ with zero average to
\begin{equation*}
-\nabla \cdot  A^{\#,N}(\xi+\nabla \phi_N^\spa)\,=\,0 \mbox{ in }\Z^d,
\end{equation*}
where $A^{\#,N}$ is the $Q_N$-periodic extension of $A|_{Q_N}$ on $\Z^d$, that is for all $k\in \Z^d$ and $z\in Q_N$, $A^{\#,N}(kN+z):=A(z)$.
The associated approximation $A_N^{\spa}$ of $A_\ho$ is then given by
\begin{equation*}
\xi\cdot A_{N}^{\spa}\xi\,:=\,\fint_{Q_N}(\xi+\nabla \phi_N^\spa) \cdot A^{\#,N}(\xi+\nabla \phi_N^\spa).
\end{equation*}
As a direct consequence of homogenization (see also \cite{Owhadi-03}), we have almost surely
\begin{equation}\label{eq:as-spa}
\lim_{N\to \infty}|A_N^\spa-A_\ho|\,=\,0.
\end{equation}

\medskip

In order to define the periodization in law, we have to make even more specific the structure of $A$.
For all $i\in \{1,\dots,d\}$, let $g_i:[\alpha,\beta]^\B\to [\alpha,\beta]$ be a measurable function depending only on a finite number of variables in $\B$ (recall that $\B$
is the set of edges), and let $\{\bar{\omega}_e\}_{e\in \B}$ be a family of i.i.d.\ random variables in $[\alpha,\beta]$.
We assume that the conductances $\{\omega_e\}_{e\in \B}$ are given as follows: For all $z\in \Z^d$ and $i\in \{1,\dots,d\}$,
$$
\omega_{(z,z+\ee_i)}\,:=\,g_i(\tau_z \bar{\omega}),
$$
where $\tau_z\bar{\omega}$ is the translation of $\bar{\omega}$ by $z$, i.~e. for all $e=(z',z'+\ee_j)$ with $z'\in \Z^d$ and $j\in \{1,\dots,d\}$, $\tau_z \bar{\omega}_e:=\bar{\omega}_{(z+z',z+z'+\ee_j)}$.
Since $g_i$ only depends on a finite number of variables, $A$ has finite correlation length.
The periodization in law consists in periodizing the underlying i.i.d.\ random variables and then applying the deterministic function $g$.
For all $N\geq 1$ we define $\bar{\omega}^{\#,N}$ by: for all $k\in \Z^d$, $z\in Q_N$, and $i\in \{1,\dots,d\}$, $\bar{\omega}^{\#,N}_{(k+z,k+z+\ee_i)}:=\bar{\omega}_{(z,z+\ee_i)}$,
and we set for all $z\in \Z^d$
$$
A_{\#,N}(z)\,:=\,\diag (g_1(\tau_z \bar{\omega}^{\#,N}),\dots,g_d(\tau_z \bar{\omega}^{\#,N})).
$$
Note that $A_{\#,N}$ is $Q_N$-periodic.
It does coincide with $A^{\#,N}$ if $g_i(\bar \omega)=\mathcal G(\bar \omega_{(0,0+\ee_i)})$ for some $\mathcal G$, in which case the conductances $\omega_e$ are i.i.d.\ random variables, but not otherwise.
We then consider the unique $Q_N$-periodic solution $\phi_N^\law$  with zero average to 
\begin{equation*}
-\nabla \cdot A_{\#,N}(\xi+\nabla \phi_N^\law)\,=\,0 \mbox{ in } \Z^d,
\end{equation*}
and define
\begin{equation*}
\xi\cdot A_{N}^{\law}\xi\,:=\,\fint_{Q_N}(\xi+\nabla \phi_N^\law) \cdot A_{\#,N}(\xi+\nabla \phi_N^\law).
\end{equation*}
In order to prove the almost sure convergence
$$
\lim_{N\to \infty} |A_N^\law-A_\ho|\,=\,0,
$$
in view of \eqref{eq:as-spa}, it is enough to show that for large $N$,
$$
\xi\cdot A_N^\spa\xi-o(1) \,\leq\, \xi\cdot A_N^\law \xi\,\leq\, \xi\cdot A_N^\spa\xi+o(1).
$$
By Meyers' estimates, there exists $p>2$ such that for all $N$ and almost surely,
\begin{equation}\label{eq:meyers}
\int_{Q_N} |\nabla \phi_N^\law|^p ,\int_{Q_N}|\nabla \phi_N^\spa|^p \,\lesssim  \, N^d.
\end{equation}
Next we use the following alternative formula for $A^\law_N$ (which holds by symmetry of $A_{\#,N}$):
$$
\xi\cdot A^\law_N \xi \,=\, \inf \left\{ \fint_{Q_N}(\xi+\nabla \phi)\cdot A_{\#,N}(\xi+\nabla \phi),\phi \text{ is }Q_N\text{-periodic}  \right\}.
$$
Using $\phi_N^\spa$ as a test function and the fact that $A_{\#,N}$ and $A^{\#,N}$ coincide on $Q_{N-\mathscr{L}_c}$, where $\mathscr{L}_c$ is the correlation-length of $A$,
we get by H\"older's inequality with exponents $(p/2,p/(p-2))$
\begin{eqnarray*}
\xi\cdot A^\law_N \xi &\leq & \fint_{Q_N}(\xi+\nabla \phi_N^\spa)\cdot A_{\#,N}(\xi+\nabla \phi_N^\spa)\\
&=& \fint_{Q_N}(\xi+\nabla \phi_N^\spa)\cdot A^{\#,N}(\xi+\nabla \phi_N^\spa)\\
&& \qquad +\fint_{Q_N}(\xi+\nabla \phi_N^\spa)\cdot (A_{\#,N}-A^{\#,N})(\xi+\nabla \phi_N^\spa)\\
&\leq & \xi\cdot A_{\#,N}\xi + \beta N^{-d} \int_{Q_N\setminus Q_{N-\calL_c}} (1+|\nabla \phi_N^\spa|^2)\\
&\leq & \xi\cdot A_{N}^\spa\xi + \beta N^{-d} |Q_N\setminus Q_{N-\calL_c}|^{(p-2)/p} \|\nabla \phi_N^\spa\|_{L^p(Q_N)}^{2/p} \\
&=&\xi\cdot A_{N}^\spa\xi +o(1)
\end{eqnarray*}
using in addition \eqref{eq:meyers}. The converse inequality is proved the same way.

\medskip

Let us turn to convergence rates, and begin with the periodization in law, which is analyzed in \cite{Gloria-Neukamm-Otto-a,Gloria-Neukamm-Otto-b}.
The general approach takes the form:
\begin{equation*}
\expec{|\xi\cdot A_N^\law\xi-\xi\cdot A_\ho\xi|^2}\,= \,\var{\xi\cdot A_N^\law\xi}+(\xi\cdot(\expec{A^\law_N}-A_\ho)\xi)^2.
\end{equation*}
As shown in \cite{Gloria-Neukamm-Otto-a,Gloria-Neukamm-Otto-b}, the random and systematic errors satisfy
\begin{equation}\label{eq:per1-random}
\var{\xi \cdot A_N^\law\xi}^{1/2}\,\lesssim \, N^{-d/2},
\end{equation}
and 
\begin{equation}\label{eq:per1-syst}
|\xi \cdot \expec{A_N^\law}\xi-\xi\cdot A_\ho\xi|\,\lesssim \, 
 N^{-d}\ln^d N.
\end{equation}
The proof of \eqref{eq:per1-syst} is subtle.
Let us give a hint of the proof in dimension $d=2,3$ for the case i.i.d.\ conductances.
The idea consists in introducing a zero-order term as for the regularization, and we consider 
\begin{equation*}
\xi\cdot A_{\mu,N}^{\law}\xi\,:=\,\fint_{Q_N}(\xi+\nabla \phi_{\mu,N}^\law) \cdot A_{\#,N}(\xi+\nabla \phi_{\mu,N}^\law),
\end{equation*}
where $\mu>0$, and $\phi_{\mu,N}^\law$ is the $Q_N$-periodic solution to 
\begin{equation*}
\mu \phi_{\mu,N}^\law-\nabla \cdot A_{\#,N}(\xi+\nabla \phi_N^\law)\,=\,0 \mbox{ in } \Z^d.
\end{equation*}
We also denote by $A_\mu$ the approximation of $A_\ho$ defined in \eqref{eq:Amu}.
Then by the triangle inequality, for all $\mu>0$,
\begin{multline*}
|\xi \cdot \expec{A_N^\law}\xi-\xi\cdot A_\ho\xi|\,\leq \,|\xi \cdot \expec{A_N^\law}\xi-\xi\cdot  \expec{A_{\mu,N}^\law}\xi|\\
+ |\xi\cdot  \expec{A_{\mu,N}^\law}\xi-\xi\cdot A_\mu \xi|+|\xi\cdot(A_\mu-A_\ho)|.
\end{multline*}
The first of the r.~h.~s.~has the same scaling as the last term (repeating the arguments of spectral theory), that is \eqref{eq:0dir-syst}, which is independent of $N$.
The only term which relates $N$ to $\mu$ is the second term. In the i.i.d.\ case, we have by periodicity and stationarity
\begin{eqnarray}
\xi\cdot \expec{A_{\mu,N}^\law}\xi&=&\expec{\fint_{Q_N}(\xi+\nabla \phi_{\mu,N}^\law) \cdot A_{\#,N}(\xi+\nabla \phi_{\mu,N}^\law)}\nonumber\\
&=&\expec{(\xi+\nabla \phi_{\mu,N}^\law(0)) \cdot A_{\#,N}^\law(0)(\xi+\nabla \phi_{\mu,N}^\law(0))}.\label{eq:coupling}
\end{eqnarray}
Hence, since $A_{\#,N}^\law(0)=A(0)$,
\begin{multline*}
\xi\cdot \expec{A_{\mu,N}^\law}\xi-\xi\cdot A_{\mu}\xi \,\\ =\, \expec{(\xi+\nabla \phi_{\mu,N}^\law(0)) \cdot A(0)(\xi+\nabla \phi_{\mu,N}^\law(0))-(\xi+\nabla \phi_{\mu}(0)) \cdot A(0)(\xi+\nabla \phi_{\mu}(0))},
\end{multline*}
and we only have to compare $\nabla \phi_{\mu}$ and $\nabla \phi_{\mu,N}^\law$ at the origin.
Using deterministic estimates on the Green's function, this yields
$$
|\xi\cdot \expec{A_{\mu,N}^\law}\xi-\xi\cdot A_{\mu}\xi |\,\lesssim \, \frac{1}{\sqrt{\mu}}\exp(-c\sqrt{\mu}N),
$$
for some $c$ depending only on $\alpha$, $\beta$, and $d$.
Optimizing the error with respect to $\mu$ (taking for instance $\sqrt{\mu} =\frac{d+1}{c}  L^{-1}\ln(L/\ln L)$) yields the announced result \eqref{eq:per1-syst}.

Let us make two comments on this proof.
First, for $d>3$ this strategy does not allow one to get \eqref{eq:per1-syst} since \eqref{eq:0dir-syst} saturates at $d=4$.
Instead of comparing $A_\ho$ to $A_\mu$, we then compare $A_\ho$ to a family of approximations $A_{k,\mu}$ for which we have
$$
\xi\cdot A_{k,\mu}\xi-\xi\cdot A_\ho\xi\,\lesssim \,\mu^p \int_{\R^+} \frac{1}{\lambda^{p+1}}de_{\mathfrak{d}}(\lambda),
$$
for some $p\in \N$ which tends to infinity as $k$ tends to infinity. This allows one to fully exploit \eqref{eq:optima-sp-exp}.
There are many possible choices for $A_{k,\mu}$, as the one introduced in \cite{Gloria-Mourrat-10a} and defined by its spectral formula (and then translated back in physical space), and the one introduced in \cite{Gloria-Neukamm-Otto-a} and defined by Richardson extrapolation in space.
Second, the proof presented above crucially relies on the i.i.d.\ structure when writing \eqref{eq:coupling}.
Indeed this relies on the fact that both statistical ensembles we use (the i.i.d.\ on $\B$ and on the set of edges of the $N$-torus) can be coupled (the ensemble on the torus being ``included" in the ensemble on $\B$). A similar construction can be made in the case of hidden i.i.d.\ variables as introduced above --- but not in full generality. 

\medskip

We expect these convergence rates to be optimal in general, so that the systematic error would scale as the square of the random error in any dimension (up to logarithmic corrections).
As for the method with the zero-order term, we expect $N^{d/2}(\xi\cdot A^\law_N\xi-\xi\cdot \expec{A^\law_N}\xi)$
to converge in law to a Gaussian random variable, this time for all $d\geq 2$.

\medskip

In the case of the periodization in space, we still expect the random error to scale as the central limit theorem
\begin{equation}\label{eq:per2-random}
\var{\xi \cdot A_N^\spa\xi}^{1/2}\,\lesssim \, N^{-d/2}.
\end{equation}
Yet, we do not expect the scaling of \eqref{eq:per1-syst} to hold in this case, and we rather conjecture that (unless the entries of $A$ are i.i.d.)\
the systematic error scales as a surface effect (as for Dirichlet boundary conditions), namely
\begin{equation}\label{eq:per2-syst}
|\xi \cdot \expec{A_{N}^\spa}\xi-\xi\cdot A_\ho\xi|\,\sim \, 
 N^{-1},
\end{equation}
although we do not have a proof of this.
The intuition behing this conjecture is that the imposed periodicity is not compatible with the underlying stationarity. A similar phenomenon occurs when considering a periodic problem and approximating the corrector on a domain which is not a multiple of the period and with periodic boundary conditions (so that the periodicity of the approximated corrector and that of the true corrector do not match), which yields an error which scales again as a surface effect (although the prefactor is usually much smaller than for Dirichlet boundary conditions).

\subsection{Numerical study}

\subsubsection{Homogeneous Dirichlet boundary conditions}
\label{ss:num1}
We consider the simplest possible example: the case of Bernoulli variables.
The conductances $\{\omega_e\}_{e\in \B}$ are i.i.d.\ random variables taking values $\alpha=1$ and
$\beta=4$ with probability $1/2$. 

\medskip

\begin{figure}
\centering
\psfrag{hh}{$\rm{log}_{10}(N)$}
\psfrag{kk}{$\rm{log}_{10}$(Syst. and rand. errors)}
\includegraphics[trim=7cm 0.5cm 7cm 0.5cm,scale=0.5]{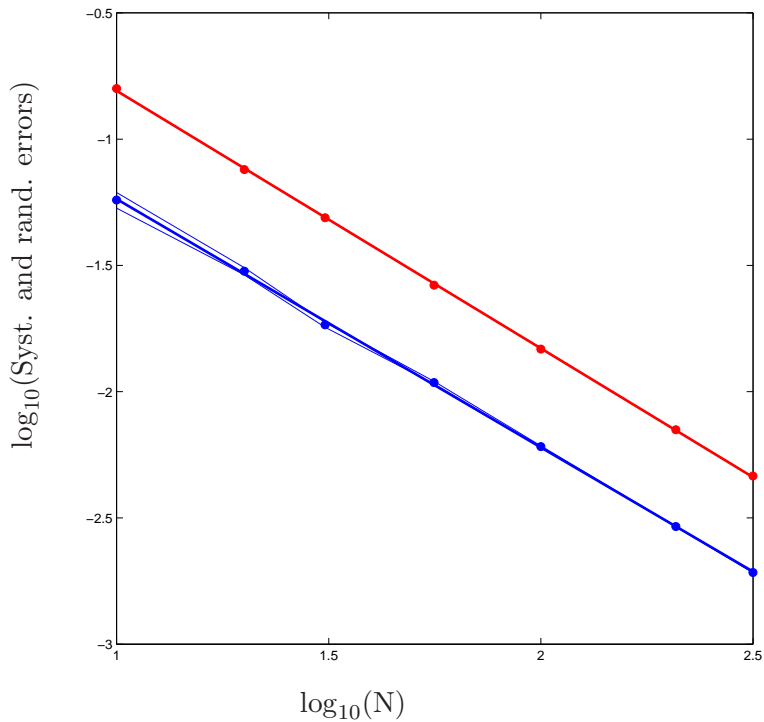}
\caption{Dirichlet boundary conditions, $d=2$,
statistical error (red) rate $1.02$ and prefactor $1.62$,  systematic error (blue) rate $0.98$ and prefactor $0.56$.}
\label{fig:2d_Dir}
\end{figure}

In dimension $d=2$, it is known that the homogenized matrix $A_\ho$ takes the form $\sqrt{\alpha \beta}\Id=2\Id$ (the Dykhne formula,
see \cite{Gloria-10} for a proof).
On Figure~\ref{fig:2d_Dir} we have plotted the estimates of the random error $N\mapsto \var{\xi\cdot A_N\xi}^{1/2}$
and of the systematic error $N\mapsto |\xi\cdot \expec{A_N}\xi-\xi\cdot A_\ho\xi|$ in logarithmic scale.
\begin{table}
\center
\caption{  \label{tab:2d_Dir}}
\begin{tabular}{|m{2cm}|c|c|c|c|c|c|c|} 
\hline
  N                       &  10   &  20   &  31   &  56  &   100  &   208  &   316    \\
\hline
  Number of realizations   &  1500  &  6000  &  14415  &  47040 &  150000  &  648960  &  1497840  \\ 
\hline
\end{tabular} 
\end{table}

These errors are approximated by empirical averages of independent realizations, and intervals of confidence are given for the systematic error (corresponding to the empirical standard deviation).
The number of independent realizations in function of $N$ is displayed for completeness in Table~\ref{tab:2d_Dir}.
As can be seen, the apparent convergence rates of the random error and of the systematic error are $1.02$ and $0.98$, respectively.
This confirms the conjecture of \eqref{eq:dir-random} and \eqref{eq:dir-syst} for $d=2$.

\medskip

\begin{figure}
\centering
\psfrag{hh}{$\rm{log}_{10}(N)$}
\psfrag{kk}{$\rm{log}_{10}$(Syst. and rand. errors)}
\includegraphics[trim=7cm 0.5cm 7cm 0.5cm,scale=0.5]{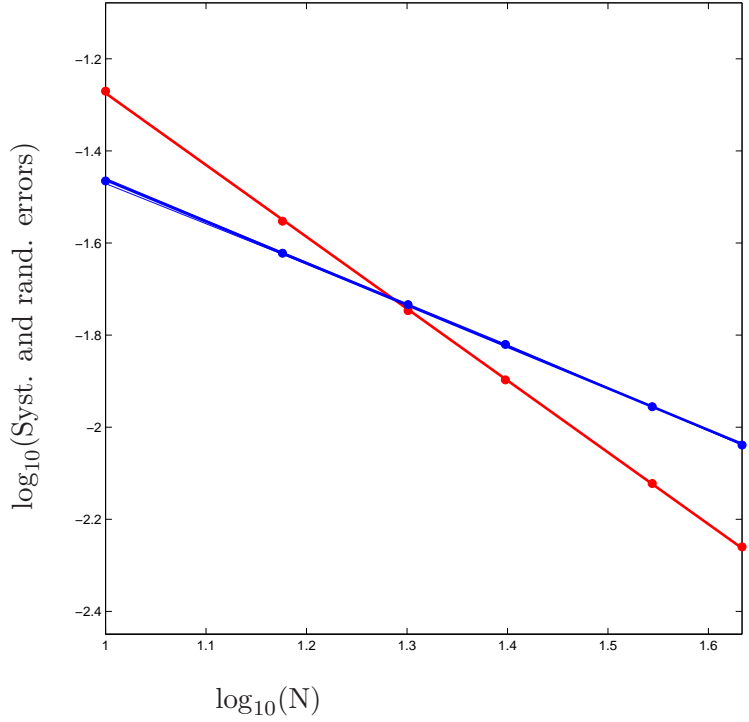}
\caption{Dirichlet boundary conditions, $d=3$,
statistical error (red) rate $1.56$ and prefactor $1.91$,  systematic error (blue) rate $0.90$ and prefactor $0.27$.}
\label{fig:3d_Dir}
\end{figure}

For numerical tests in dimension $d=3$, we have to proceed slightly differently since there is no closed formula for the homogenized coefficient (the homogenized matrix  is still a multiple of the identity by symmetry arguments). 
The approximation of the random error $N\mapsto \var{\xi\cdot A_N\xi}^{1/2}$ is unchanged.
Yet, instead of plotting the systematic error $N\mapsto  |\xi\cdot \expec{A_N}\xi-\xi\cdot A_\ho\xi|$, we plot an approximation 
of $N\mapsto  |\xi\cdot \expec{A_N}\xi-\xi\cdot \expec{A_N^\law}\xi|$ via empirical averages, where $A_N^\law$ is the approximation of $A_\ho$ by periodization.
In particular, in view of \eqref{eq:per1-syst}, we have by the triangle inequality
$$
|\xi\cdot \expec{A_N}\xi-\xi\cdot A_\ho\xi| \,\geq \, |\xi\cdot \expec{A_N}\xi-\xi\cdot \expec{A_N^\law}\xi|-CN^{-d}\ln^dN
$$
for some $C>0$, so that $|\xi\cdot \expec{A_N}\xi-\xi\cdot A_\ho\xi| $ and $ |\xi\cdot \expec{A_N}\xi-\xi\cdot \expec{A_N^\law}\xi|$
are of the same order provided $ |\xi\cdot \expec{A_N}\xi-\xi\cdot \expec{A_N^\law}\xi|\gtrsim N^{-d}\ln^dN$ (which we indeed observe numerically).
These two errors $N\mapsto \var{\xi\cdot A_N\xi}^{1/2}$ and $N\mapsto  |\xi\cdot \expec{A_N}\xi-\xi\cdot \expec{A_N^\law}\xi|$ are plotted on 
Figure~\ref{fig:3d_Dir} in logarithmic scale.
The number of independent realizations in function of $N$ is displayed for completeness in Table~\ref{tab:3d_Dir}.
\begin{table}
\center
\caption{  \label{tab:3d_Dir}}
\begin{tabular}{|m{2cm}|c|c|c|c|c|c|} 
\hline
  N                       &  10   &  15   &  20   &  25  &      35  &   43    \\
\hline
  Number of realizations   &  15000  &  50625  &  120000  &  234375 &  643125  & 1192605   \\ 
\hline
\end{tabular} 
\end{table}
As can be seen, the apparent convergence rate of the random error and of the (modified) systematic error are $1.56$ and $0.9$, respectively. 
This confirms the conjecture of \eqref{eq:dir-random} and \eqref{eq:dir-syst} for $d=3$.

%

\subsubsection{Regularized corrector and filtering}\label{sec:regu}

\begin{figure}
\centering
\psfrag{hh}{$\rm{log}_{10}(N)$}
\psfrag{kk}{$\rm{log}_{10}$(Syst. and rand. errors)}
\includegraphics[trim=7cm 0.5cm 7cm 0cm,scale=0.5]{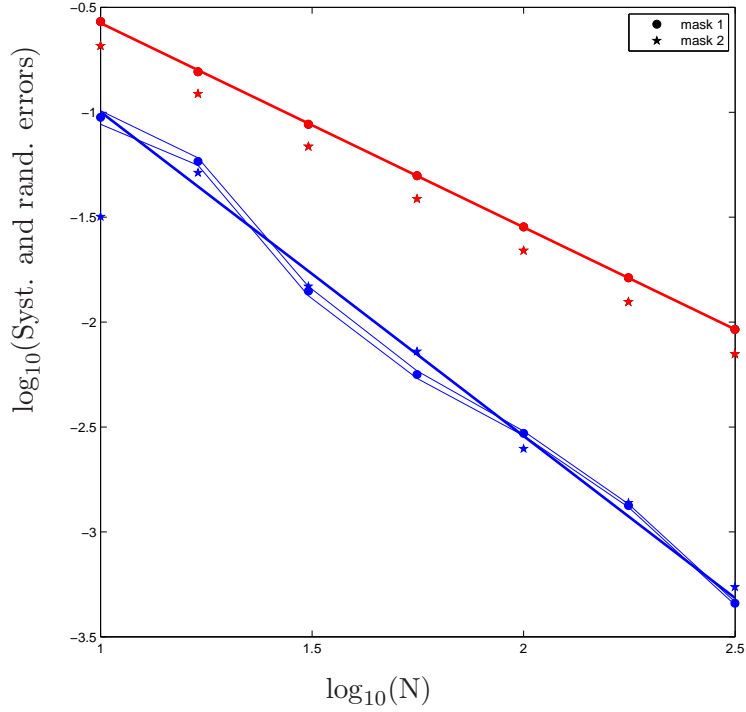}
\caption{Regularized corrector, $d=2$, i~.i.~d. case,
statistical error (red) rate $0.97$ and prefactor $2.51$,  systematic error (blue) rate $1.55$ and prefactor $3.55$.}
\label{fig:2d_zero_stat}
\end{figure}

\begin{figure}
\centering
\includegraphics[trim=9cm 0.5cm 7cm 0.5cm,scale=0.4]{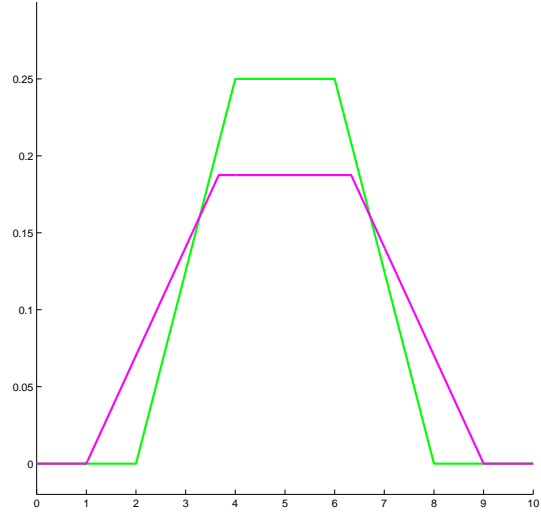}
\caption{Two masks for $N=10$}
\label{fig:mask}
\end{figure}

We still consider the same simple two-dimensional example of the Bernoulli random variables taking values $\alpha=1$ and
$\beta=4$ with probability $1/2$, and for which the homogenized matrix is $A_\ho=\sqrt{\alpha\beta}\Id=2\Id$.
In order to define $A_{\mu,N,L}$ completely, we need to choose $L$ and $\mu$ in function of $N$, and define the averaging mask $\eta_L$.
We have taken
\begin{itemize}
\item $L=4N/5$ and $L=3N/5$,
\item $\mu=125/N^{3/2}$,
\item a piecewise affine mask, as plotted on Figure~\ref{fig:mask} for $N=10$ (the first one for $L=4N/5$ and the second one for $L=3N/5$).
\end{itemize}
The theoretical predictions \eqref{eq:0dir-random} and \eqref{eq:0dir-syst} take the following forms with these parameters:
\begin{eqnarray}
\var{\xi \cdot A_{\mu,N,L}\xi}^{1/2} &\lesssim & N^{-1}\ln^q N\label{eq:0dir-random-test}\\
|\xi\cdot \expec{A_{\mu,N,L}}\xi-\xi\cdot A_\ho\xi|& \lesssim &  N^{-3/2} \ln^2 N.\label{eq:0dir-syst-test}
\end{eqnarray}
These two errors are plotted on 
Figure~\ref{fig:2d_zero_stat} in logarithmic scale.
The number of independent realizations in function of $N$ is displayed for completeness in Table~\ref{tab:2d_zero}.
\begin{table}
\center
\caption{  \label{tab:2d_zero}}
\begin{tabular}{|m{2cm}|c|c|c|c|c|c|c|} 
\hline
  N                       &  10     &  17    &  31     &  56    &   100     &   177     &  316      \\
\hline
  Number of realizations   &  1500  &  4335  &  14415  &  47040 &  150000   &  510081   & 1480338 \\ 
\hline
\end{tabular} 
\end{table}
As can be seen, the apparent convergence rates for the random and systematic errors are $0.97$ and $1.55$, respectively.
This shows the sharpness of the analysis.
These numerical tests also give an idea on the prefactors in \eqref{eq:0dir-random} and \eqref{eq:0dir-syst}. We observe that the systematic error decays faster than the random error, and that the prefactors are of the same order (roughly twice as big for the systematic error).

\medskip

\begin{figure}
\centering
\psfrag{hh}{$\rm{log}_{10}(N)$}
\psfrag{kk}{$\rm{log}_{10}$(Syst. and rand. errors)}
\includegraphics[trim=7cm 0.5cm 7cm 0cm,scale=0.5]{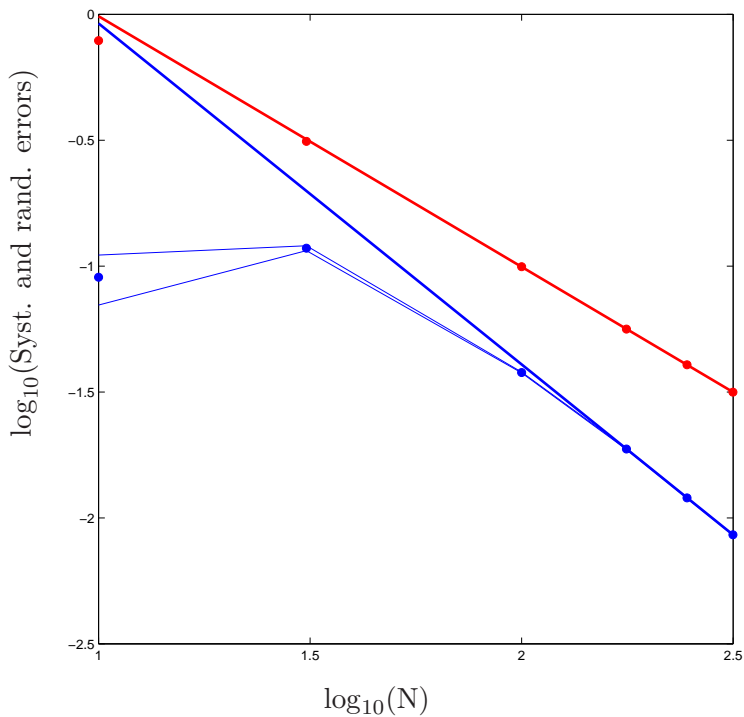}
\caption{Regularized corrector, $d=2$, correlated case,
statistical error (red) rate $1$ and prefactor $9.77$,  systematic error (blue) rate $1.35$ and prefactor $20.9$.}
\label{fig:2d_zero_corr}
\end{figure}

The second series of tests on the regularization method aims at validating an approach which will be used to compare the two periodization methods.
In particular, we consider now a two-dimensional case with correlations for which we do not have a closed formula of the homogenized coefficients.
The statistics of the coefficients is defined as follows.
We let $\{\bar \omega_{(z,z+e_i)}\}_{z\in \Z^2,i\in\{1,2\}}$ be i.i.d.\ variables following a uniform law in $[0,1]$.
We define $\omega_{(z,z+e_i)}$ to be $\alpha=1$ if for all $z'$ such that $\|z'-z\|_{\infty}\leq 2$ we have $\bar \omega_{(z,z+e_i)}\leq p$ (that is, if the 25 hidden i.i.d.\ random variables are less than $p$),
and  $\omega_{(z,z+e_i)}$ to be $\beta=4$ if there exists $z'$ with $\|z'-z\|_{\infty}\leq 2$ such that $\bar \omega_{(z,z+e_i)}>p$,
where $p$ is chosen so that $p^{25}=1/2$ (that is, $\alpha$ and $\beta$ are equiprobable).
The typical realization of such conductances is made of islands of $\beta$'s of size 4 in a sea of $\alpha$'s. 
This example can be recast in the form of the correlated case described in Paragraph~\ref{ss:period}.
Since the analysis of \cite{Gloria-Otto-09b,Gloria-Neukamm-Otto-a,Gloria-Neukamm-Otto-b} also covers the case of finite correlation length,  the theoretical predictions \eqref{eq:0dir-random-test} and \eqref{eq:0dir-syst-test}
also apply to this case, provided we still take
\begin{itemize}
\item $L=4N/5$;
\item $\mu=125/N^{3/2}$,
\item a piecewise affine mask, as plotted on Figure~\ref{fig:mask}.
\end{itemize}
Since we do not know $A_\ho$ a priori, we shall replace the systematic error by $N\mapsto |\xi\cdot \expec{A_{\mu,N,L}-A_N^\law}\xi|$, where
$A_N^\law$ is the approximation of $A_\ho$ by the periodization in law method.
The combination of \eqref{eq:0dir-syst-test} and \eqref{eq:per1-syst} indeed yields
$$
|\xi\cdot \expec{A_{\mu,N,L}-A_N^\law}\xi|\,\lesssim\, N^{-3/2}\ln^2N,
$$
which we want to check numerically.
This modified systematic error is plotted on 
Figure~\ref{fig:2d_zero_corr} in logarithmic scale.
The number of independent realizations in function of $N$ is reported on in Table~\ref{tab:2d_reg_corr}.
\begin{table}
\center
\caption{  \label{tab:2d_reg_corr}}
\begin{tabular}{|m{2cm}|c|c|c|c|c|c|
} 
\hline
  N                       &  10      &  31     
&  100      & 177   & 246 & 316
\\
\hline
  Number of realizations   &  1500  & 14415  & 150000  & 469935 & 907740  &   1497840
 \\ 
\hline
\end{tabular} 
\end{table} 
As can be seen, the apparent convergence rate for the modified systematic error is close to $3/2$, so that
the true systematic error $|\xi\cdot \expec{A_{\mu,N,L}}-A_\ho\xi|$ has the same decay (up to possible logarithmic corrections) 
since
$$
|\xi\cdot \expec{A_{\mu,N,L}}-A_\ho\xi|\,\geq \, |\xi\cdot \expec{A_{\mu,N,L}-A_N^\law}\xi|+CN^{-2}\ln^2N,
$$
for some $C>0$ due to \eqref{eq:per1-syst}.
Note that the asymptotic regime is more difficult to capture in the correlated case than in the i.i.d.\ since
the typical lengthscale is 4 in the first case (the size of a typical island), and 1 in the second case.
Yet these tests show it is possible to observe numerically a convergence with a rate larger than $1$ in this correlated case.

\subsubsection{Periodization methods}

\begin{figure}
\centering
\psfrag{hh}{$\rm{log}_{10}(N)$}
\psfrag{kk}{$\rm{log}_{10}$(Syst. and rand. errors)}
\includegraphics[trim=7cm 0.5cm 7cm 0.5cm,scale=0.5]{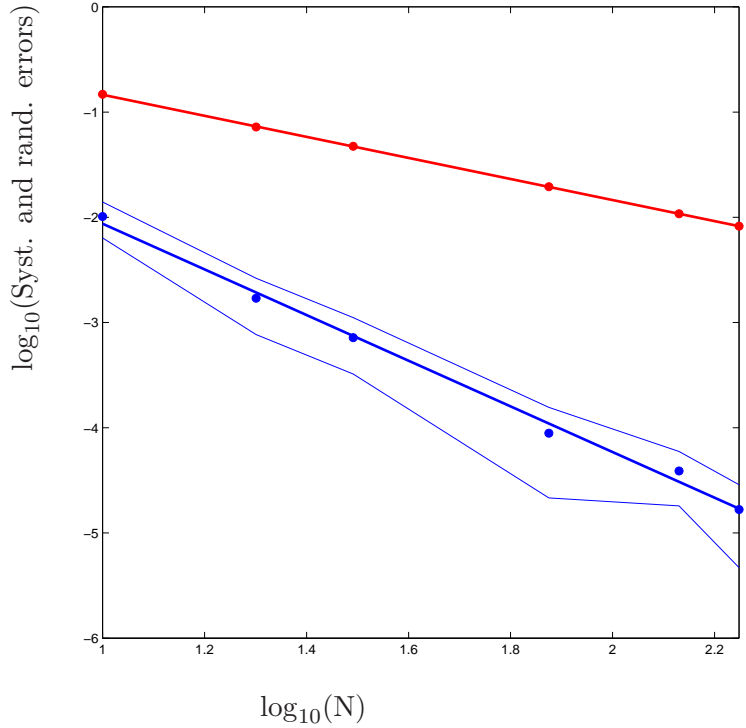}
\caption{Periodization in law, $d=2$, i.i.d.,\ 
statistical error (red) rate $1.01$ and prefactor $1.48$,  systematic error (blue) rate $2.17$ and prefactor $1.29$.}
\label{fig:2d_per}
\end{figure}

In this subsection we first check numerically the sharpness of \eqref{eq:per1-random} and \eqref{eq:per1-syst} on our simple
two-dimensional Bernoulli case.
On Figure~\ref{fig:2d_per} we have plotted the estimates of the random error $N\mapsto \var{\xi\cdot A_N^\law\xi}^{1/2}$
and of the systematic error $N\mapsto |\xi\cdot \expec{A_N^\law}\xi-\xi\cdot A_\ho\xi|$ in logarithmic scale.
\begin{table}
\center
\caption{  \label{tab:2d_per}}
\begin{tabular}{|m{2cm}|c|c|c|c|c|c|
} 
\hline
  N                       &  10     &  20    &  31     
&  75           & 135    & 177   
\\
\hline
  Number of realizations   &1500    &   14415    &  498904  &    469930   &   780900   &  1245782
%
%
 \\ 
\hline
\end{tabular} 
\end{table}
These errors are approximated by empirical averages of independent realizations, and intervals of confidence are given for the systematic error (corresponding to the empirical standard deviation).
The number of independent realizations in function of $N$ is displayed for completeness in Table~\ref{tab:2d_per}.
As can be seen, the apparent convergence rates of the random error and of the systematic error are $1.01$ and $2.17$ (note that the fluctuations are more important for the evaluation of the systematic error which is very small), which confirms the predictions (up to the logarithmic correction for the systematic error).
In addition, the prefactors are again of the same order (slightly smaller for the systematic error), so that the systematic error is really negligible in front of the random error.
This will be further analyzed in the next subsection.


\medskip

\begin{figure}
\centering
\psfrag{hh}{$\rm{log}_{10}(N)$}
\psfrag{kk}{$\rm{log}_{10}$(Syst. and rand. errors)}
\includegraphics[trim=9cm 0.5cm 7cm 0.5cm,scale=0.5]{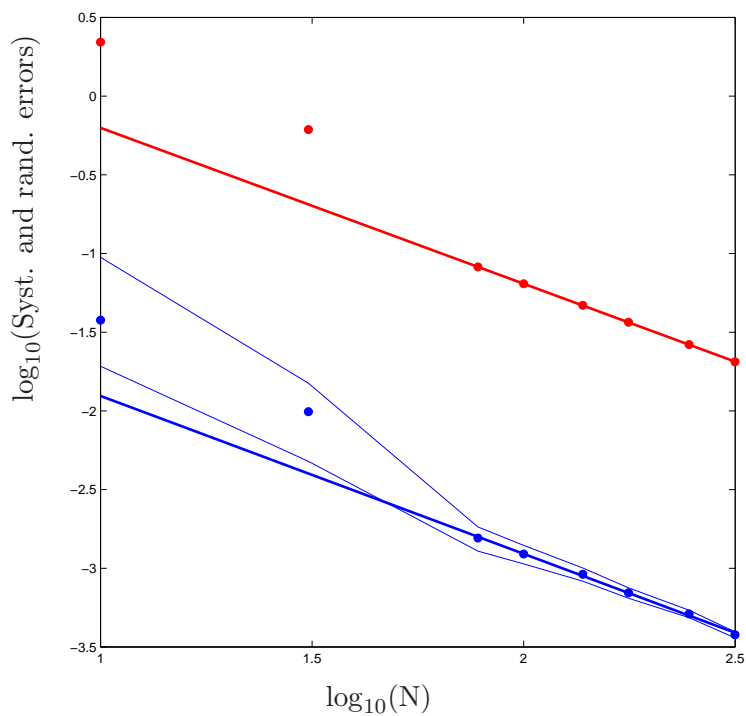}
\caption{Periodization in space, $d=2$, correlated,
statistical error (red) rate $0.99$ and prefactor $6.17$,  modified systematic error (blue) rate $1$ and prefactor $0.13$.}
\label{fig:2d_malcorr_1}
\end{figure}

The second series of tests deals with the correlated two-dimensional example introduced in Paragraph~\ref{sec:regu}.
In this case, we'd like to check numerically our conjecture \eqref{eq:per2-syst} on the systematic error for the periodization in space method.
As in the previous section, we replace the systematic error by the modified systematic error 
$N\mapsto |\xi\cdot \expec{A^\spa_N-A^\law_N}\xi|$.
In view of \eqref{eq:per1-syst}, we indeed have
$$
 |\xi\cdot (\expec{A^\spa_N}-A_\ho)\xi|\,=\, |\xi\cdot \expec{A^\spa_N-A^\law_N}\xi|+O(N^{-2}\ln^2N).
$$
On Figure~\ref{fig:2d_malcorr_1} we have plotted the estimates of the random error $N\mapsto \var{\xi\cdot A_N^\law\xi}^{1/2}$
and of the modified systematic error $N\mapsto |\xi\cdot \expec{A_N^\law -A_N^\spa}\xi|$ in logarithmic scale.
\begin{table}
\center
\caption{  \label{tab:2d_per-mal}}
\begin{tabular}{|m{2cm}|c|c|c|c|c|c|c|c|} 
\hline
  N                             &  10& 31& 78        & 100    & 138    & 177    &  246 &316   \\
\hline
  Number of realizations       & 1500 & 14415 &  91260        & 150000 &285660 & 469935 & 907740  &   1497840\\ 
\hline
\end{tabular} 
\end{table}
These errors are approximated by empirical averages of independent realizations, and intervals of confidence are given for the systematic error (corresponding to the empirical standard deviation).
The number of independent realizations in function of $N$ is displayed for completeness in Table~\ref{tab:2d_per-mal}.
As can be seen, the slopes of the random error and of the systematic error are $0.99$ and $1$, which confirms the conjectures \eqref{eq:per2-random} and \eqref{eq:per2-syst}
on the periodization in space method.
Yet it should be noticed that the prefactor of the systematic error is much smaller (around 40 times smaller) than the prefactor of the random error.
This makes the systematic error hardly noticeable in practice. This would'nt be so clear in the case of dimension $d\geq 3$ (as on Figure~\ref{fig:3d_Dir}).

\medskip

\begin{figure}
\begin{center}
\begin{minipage}{0.48\textwidth}
\includegraphics[scale=0.5]{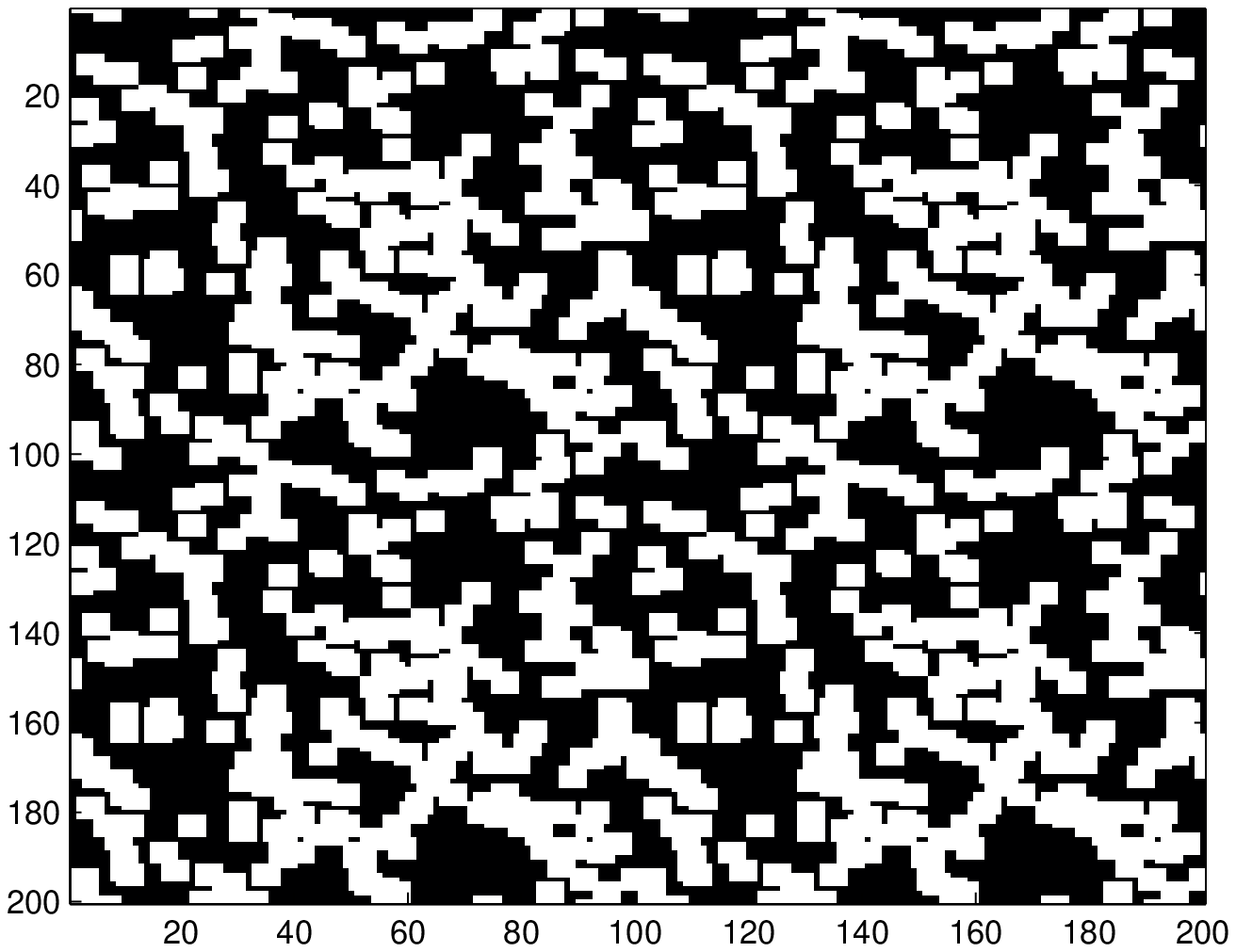}
\caption{Periodization in law of $A$ (four periods)}
\label{fig:A-good}
\end{minipage}
\begin{minipage}{0.48\textwidth}
\includegraphics[scale=0.5]{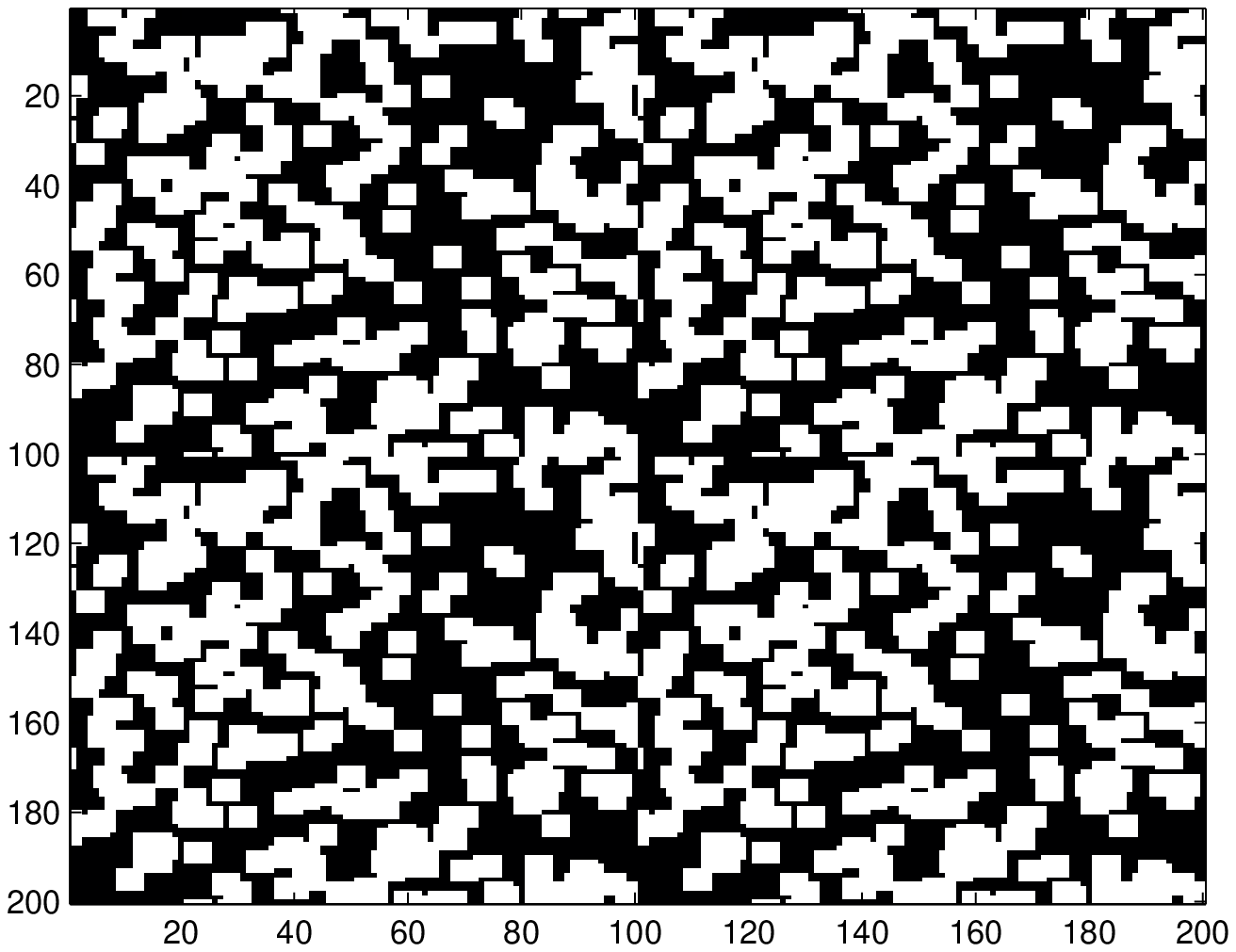}
\caption{Periodization in space of $A$ (four periods)}
\label{fig:A-bad}
\end{minipage}
\end{center}
\end{figure}

It is instructive to compare the periodization in law $A_{N,\#}$ of $A$ to the periodization in space $A^{N,\#}$ of $A$ in terms of images.
Typical realizations of these periodizations are represented on Figures~\ref{fig:A-good} and~\ref{fig:A-bad} for $N=100$ --- 4 periods are reproduced.
A close look at Figure~\ref{fig:A-bad} reveals the mismatch between stationarity and enforced periodicity (which does not appear on Figure~\ref{fig:A-good} --- this is the core of the ``coupling" strategy).

\subsubsection{Optimal numerical strategy}

In this last paragraph, we propose a numerical strategy to obtain an approximation of the homogenized coefficients in the i.i.d.\ case at a given precision and
at the lowest cost.
In view of the previous subsections, it is clear that the periodization (in law) method minimizes both the random and systematic errors at $N$ fixed.
One feature we have not used yet is the scalability of the variance: if $\{A_i\}_{1\leq i \leq k}$ are $k$ independent realizations of a random variable $A$,
then
$$
\var{\frac{1}{k} \sum_{i=1}^k A_i} \,=\, \frac{1}{\sqrt{k}} \var{A}.
$$
While empirical averages of independent realizations do not allow one to reduce the systematic error, they do allow one to reduce the random error. 
This is particularly interesting since: 
\begin{itemize}
\item the dominant error is the random error,
\item it is computationally cheaper to solve several smaller linear problems than one single large problem (the solution cost of a linear system is always superlinear).
\end{itemize}

In view of the analysis of the periodization in law, in order to approximate homogenized coefficients  within a tolerance $\delta>0$,
the cheapest computational way consists in solving $k_\delta$ independent periodic problems of size $N_\delta$, and take as approximation the empirical average of these $k_\delta$ realizations $A_{N_\delta,i}^\law$:
\begin{equation}\label{eq:best}
\frac{1}{k_\delta} \sum_{i=1}^{k_\delta}A_{N_\delta,i}^\law,
\end{equation}
where $N_\delta$ is such that
$$
\delta/2 \,=\, C_{syst} N^{-d}_\delta \ln^d N_\delta,
$$
and $k_\delta$ such that
$$
\delta/2 \,=\, \frac{1}{\sqrt{k_\delta}}C_{rand} N^{-d/2}_\delta.
$$
where $C_{syst}$ and $C_{rand}$ are the optimal prefactors in the estimates \eqref{eq:per1-syst} and \eqref{eq:per1-random}, respectively.
Then we have
$$
\expec{\left|A_\ho-\frac{1}{k_\delta} \sum_{i=1}^{k_\delta}A_{N_\delta,i}^\law\right|^2}^{1/2}\,\leq \, \delta.
$$
To be precise, the associated computational cost is $k_\delta \gamma_d(N_\delta^d)$, where $M\mapsto \gamma_d(M)$ is the cost of solving a linear problem of  order $M$ (order of the matrix of the linear system) in dimension $d$. 
Since $\gamma_d$ is superlinear, one readily convinces oneself that the best one can do is indeed \eqref{eq:best}.

In the i.i.d.\ case of Figure~\ref{fig:2d_per}, an approximation of these prefactors is
$$
C_{syst}=1.29, \quad C_{rand}=1.48.
$$

\section{Numerical approximation of the homogenized coefficients using the RWRE}
\label{s:num-RWRE}
\setcounter{equation}{0}
\subsection{General approach}

We now discuss how to approximate homogenized coefficients by simulating random walks. The simulation of a random walk has a very interesting feature: one does not need to generate a full environment a priori. Rather, it suffices to generate the environment along the trajectory of the random walk. This is particularly interesting in dimensions $3$ and higher, where the random walk is transient, and visits only a vanishing fraction of the space. In fact, although the walk is recurrent in dimension~$2$, this last observation still holds in this case, since the time necessary to exit a box of size $N$ is of order $N^2$, while at such a time, the walk has typically visited a number of distinct sites of order $N^2/\log(N)$.

In general terms, the strategy consists in simulating a large number of random walks, each in its own independent environment, and rely on Theorem~\ref{t:kv} to recover the homogenized coefficients. Keeping the environment fixed would be more difficult to analyse from a theoretical point of view, would certainly lead to larger errors (although we cannot prove this), and would force us to abandon the approach of generating the environment along the trajectory.

\subsection{The discrete-time random walk}
Although it is easier to see the link between the corrector equation and the continuous-time random walk, when it comes to simulations, it is more convenient to work with a discrete-time version of $X$, since there is no waiting times to compute. We define $(Y_n)_{n \in \N}$ to be the discrete-time Markov chain such that $\PPo_z[X_0 = z] = 1$ and
\begin{equation}
\label{def:Y}
\PPo_z[Y_1 = z'] = 
\left|
\begin{array}{ll}
p(z\leadsto z') & \text{if } z' \sim z, \\
0 & \text{otherwise},	
\end{array}
\right.
\end{equation}
where $p(z\leadsto z')$ was introduced in \eqref{def:leadsto}. Considering the ``constructive'' definition of the random variable $X$ given in paragraph~\ref{ss:ctrw}, one can see that if $Y_n$ is the position of $X$ after $n$ steps, then $(Y_n)_{n \in \N}$ is a discrete-time Markov process satisfying \eqref{def:Y}. As before, it will be convenient to consider the environment viewed by the walk $Y$, defined as
$$
\omega_n = \theta_{Y_n} \ \om.
$$
Note that the probabilities $p(z\leadsto z')$ and $p(z'\leadsto z)$ need not be equal, hence it is not true in general that the counting measure is reversible for $Y$. A reversible measure $\pi$ for $Y$ should satisfy, for any $z,z' \in \Z^d$,
\begin{equation*}
\pi(z) \ p(z\leadsto z') = \pi(z') \ p(z'\leadsto z).
\end{equation*}
This relation holds if we choose $\pi(z)=p_\omega(z)$ (recall the definitions of $p_\omega(z)$ and $p(z\leadsto z')$ given in \eqref{defpomega} and \eqref{def:leadsto} respectively). This means that contrary to the continuous-time random walk, the discrete-time walk will preferably spend time on sites where $p_\omega(z)$ is large. This effect has its counterpart concerning the environment viewed by $Y$: the measure $\P$ is not invariant for $(\omega_n)_{n \in \N}$ in general, but rather the following ``tilted'' version $\td{\P}$, defined by
\begin{equation}
\label{deftdP}
\d \td{\P}(\omega) = \frac{p(\omega)}{\expec{p}} \ \d \P(\omega),
\end{equation}
where for simplicity we write $p(\omega) = p_\omega(0)$, and $\expec{p} = \expec{p(\omega)}$. In particular, the measure~$\td{\P}$ has a density with respect to our initial measure $\P$. Since $p(\omega)$ is never equal to $0$ or infinity, a property holds $\td{\P}$-almost surely if and only if it holds with $\P$-almost surely. 

We write $\td{\P}_0$ for the measure $\td{\P}\PPo_0$, and $\td{\E}_0$ for the associated expectation. Adapting slightly the arguments of paragraph~\ref{ss:clt}, we obtain the following.
\begin{thm}[\cite{Kipnis-Varadhan-86}]
\label{t:kvdisc}
Under the measure $\td{\P}_0$ and as $\eps$ tends to $0$, the rescaled discrete-time random walk $Y^{(\eps)} := (\sqrt{\e}Y_{\lfloor t/\e\rfloor})_{t \in \R_+}$ converges in distribution (for the Skorokhod topology) to a Brownian motion with covariance matrix $2\Ahd$, where 
\begin{equation}
\label{defAhd}
\Ahd = \expec{p}^{-1}  \Ah = (2d\expec{\om_e})^{-1}  \Ah,
\end{equation}
and $\Ah$ is as in \eqref{2-1:eq:hom-coeff}. In other words, for any bounded continuous functional $F$ on the space of cadlag functions, one has
\begin{equation}
\label{convergegendisc}
\td{\E}_0\Ll[ F(Y^{(\eps)}) \Rr] \xrightarrow[\eps \to 0]{} E[F(B)],
\end{equation}
where $B$ is a Brownian motion started at the origin and with covariance matrix $2 \Ahd$.  Moreover, for any $\xi \in \R^d$, one has
\begin{equation}
\label{convergesquaredisc}
n^{-1} \td{\E}_0\Ll[ \Ll(\xi \cdot Y_n\Rr)^2 \Rr] \xrightarrow[n \to + \infty]{} 2 \xi \cdot \Ahd \xi.
\end{equation}
\end{thm}
Before discussing numerical methods based on this theorem, we introduce some notation. Let $Y^{(1)}, Y^{(2)}, \ldots$ be independent discrete-time random walks evolving in the environments $\omega^{(1)}, \omega^{(2)}, \ldots$ respectively. We write $\PP^{\ov{\om}}_0$ for their joint distribution, all random walks starting from $0$, where $\ov{\om} = (\om^{(1)},\om^{(2)},\ldots)$. The environments $(\omega^{(1)}, \omega^{(2)}, \ldots)$ are themselves i.i.d.\ random variables distributed according to $\P$, and we write $\P^\otimes$ for their joint distribution. We also write $\P^\otimes_0$ as a shorthand for the measure $\P^\otimes\PP^{\ov{\om}}_0$. As usual, we simply replace ``P'' by ``E'' with the appropriate typography to denote corresponding expectations.

\subsection{Methods and theoretical analysis}
In the two following paragraphs, we discuss numerical methods based respectively on the convergences in \eqref{convergesquaredisc} and \eqref{convergegendisc}. We do so with the aim of explaining the quantitative estimates on the errors that can be proved, and the theoretical reasons lying behind the superiority of the method based on the convergence of the rescaled mean square displacement.

\subsubsection{Method based on the mean square displacement}
As announced, we start with a method based on the convergence in \eqref{convergesquaredisc}. Recall that, by the definition of the tilted measure $\td{\P}$ in \eqref{deftdP}, we have
\begin{equation}
\label{eqtilt}
n^{-1} \td{\E}_0\Ll[ \Ll(\xi \cdot Y_n\Rr)^2 \Rr] = \frac{1}{n \expec{p}} \E_0\Ll[p(\omega) (\xi \cdot Y_n)^2\Rr].
\end{equation}
By the law of large numbers, we know that for any fixed $n$, the quantity
\begin{equation}
\label{defhatA}
\hat{A}_k(n) := \frac{p(\omega^{(1)}) (\xi \cdot Y_n^{(1)})^2 + \cdots + p(\omega^{(k)}) (\xi \cdot Y_n^{(k)})^2}{kn \expec{p}} 
\end{equation}
converges (almost surely) to the quantity in the right-hand side of equality \eqref{eqtilt} as $k$ tends to infinity. From the convergence in \eqref{convergesquaredisc}, we thus obtain
\begin{equation}
\label{asymptothatA}
\lim_{n \to +\infty} \lim_{k \to + \infty} \hat{A}_k(n) = 2 \xi \cdot \Ahd \xi.
\end{equation}
The quantity $\hat{A}_k(n)$ is what we will compute. It involves $k$ random walks, each simulated in its own environment, up to time $n$. The expression involves the average $\expec{p} = 2d \expec{\omega_e}$, but this can be easily computed beforehand, so we assume that we have exact knowledge of this quantity. The convergence in \eqref{asymptothatA} can be made quantitative. 
\begin{thm}[\cite{Gloria-Mourrat-10b}]
\label{t:MCdisc}
There exist constants $q,C,c > 0$ such that for any $k \in \N^*$, any $\eps > 0$ and any $n$ large enough,
\begin{equation}\label{eq:complete}
\P^\otimes_0\left[\big|\hat{A}_k(n) - 2 \xi \cdot \Ahd \xi \big| \ge (C\mu_d(n)+\eps) /n\right] \le  \exp\left( - \frac{k \eps^2}{c n^2}\right),
\end{equation}
where
$$
\mu_d(n)=\left|
\begin{array}{ll} 
\ln^q n &  \text{if } d = 2, \\ 
1 & \text{if } d \ge 3.
\end{array}\right.
$$ 
\end{thm}
This theorem ensures that it suffices to choose $k$ as a large constant times $n^2$ to guarantee, with probability close to $1$, that the difference between $\hat{A}_k(n)$ and $2 \xi \cdot \Ahd \xi$ is smaller than $C/n$ (or $C \log^q(n)/n$ if $d = 2$) for some constant $C$.

We decompose the error
$$
\big|\hat{A}_k(n) - 2 \xi \cdot \Ahd \xi \big|
$$
into the sum of
\begin{equation}
\label{decomperrorin2}
\Ll|\hat{A}_k(n) - n^{-1} \td{\E}_0\Ll[ \Ll(\xi \cdot Y_n\Rr)^2 \Rr] \Rr| + \Ll| n^{-1} \td{\E}_0\Ll[ \Ll(\xi \cdot Y_n\Rr)^2 \Rr] - 2\xi \cdot \Ahd \xi \Rr|.
\end{equation}
We call the first difference the \emph{statistical error}, and the second difference the \emph{systematic error}. The theorem is proved if we show the following two inequalities:
\begin{equation}
\label{staterr}
\P_0^\otimes\Ll[ \Ll|\hat{A}_k(n) - n^{-1} \td{\E}_0\Ll[ \Ll(\xi \cdot Y_n\Rr)^2 \Rr] \Rr| \ge \eps/n \Rr] \le \exp\left( - \frac{k \eps^2}{c n^2}\right),
\end{equation}
\begin{equation}
\label{systerr}
\Ll| n^{-1} \td{\E}_0\Ll[ \Ll(\xi \cdot Y_n\Rr)^2 \Rr] - 2\xi \cdot \Ahd \xi \Rr| \le C\mu_d(n)/n.
\end{equation}
Inequality~\eqref{staterr} is controlled using classical large deviations theory, noting that $\hat{A}_k(n)$ is a sum of i.i.d.\ random variables. In the same vein, it is also possible to show (see \cite[Proposition~5.1]{Gloria-Mourrat-10b}) that for any sequence $k_n$ tending to infinity with $n$,
\begin{equation}
\label{convgauss}
\sqrt{k_n} \Ll( \hat{A}_{k_n}(n) - n^{-1} \td{\E}_0\Ll[ \Ll(\xi \cdot Y_n\Rr)^2 \Rr] \Rr) \xrightarrow[n \to + \infty]{\text{(distr.)}} \mcl{N}(0,\mathsf{v}),
\end{equation}
where $\mcl{N}(0,\mathsf{v})$ is a Gaussian random variable of variance $\mathsf{v}$, and
\begin{equation}
\label{defmathsfv}
\mathsf{v} = \left( 3 \frac{\expec{p^2}}{\expec{p}^2} - 1 \right) (2 \xi \cdot \Ahd \xi)^2.
\end{equation}
Although the left-hand side of inequality \eqref{systerr} is deterministic, its proof requires a more subtle analysis. The problem can be rephrased by saying that a quantitative control of the convergence in \eqref{convergesquaredisc} is needed.

For simplicity, we will instead discuss how to obtain a quantitative control of the convergence in \eqref{convergesquare}, that is, for the continuous-time random walk. 
Recall that in paragraph~\ref{ss:clt}, we decomposed $\xi \cdot X_t$ into $M_t + R_t$, where $M_t$ is a stationary martingale under the measure $\ov{\P}_0$, and $R_t$ is a remainder. The martingale property and the stationarity of the increments guarantee that $\E_0[M_t^2]$ grows linearly with $t$, and in fact (see \eqref{meanlinear} and \eqref{identifcv})
\begin{equation}
\label{miracle1}
\E_0[M_t^2] = 2 t \ \xi \cdot \Ah \xi.
\end{equation}
We stress that the fact that there is an equality above, and not just an asymptotic equivalence, is in itself a ``miracle'': it says that for any $t$, the expectatation of $M_t^2$ is \emph{exactly} the expectation of $(\xi \cdot B_t)$, where $B$ is the limiting Brownian motion (and moreover, note that we did not need any quantitative information whatsoever when deriving this identity).

Using \eqref{miracle1} and the fact that $\xi \cdot X_t = M_t+R_t$, we thus obtain
\begin{equation}
\label{decompsquare}
\E_0\Ll[(\xi \cdot X_t)^2\Rr] = 2t \xi \cdot \Ah \xi + \E_0\Ll[R_t^2\Rr] + 2 \E_0\Ll[ M_t R_t \Rr].
\end{equation}
In paragraph~\ref{ss:clt}, we have sketched the argument leading to a proof of the fact that $R_t/\sqrt{t}$ tends to $0$ in $L^2(\P_0)$ (see \eqref{2-1:eq:KV-2}): the problem is reduced to a proof of the fact that the spectral measure is light enough in the neighborhood of $0$ to integrate $\lambda^{-1}$ (see \eqref{2-1:eq:KV-1}); and this in turn is proved via a general argument of symmetry.

This general argument is interesting, and covers arbitrary stationary distributions of conductances. However, it gives no information on the rate at which $t^{-1} \E_0\Ll[R_t^2\Rr]$ converges to $0$, and this is inherent to such a high level of generality. When the conductances are independent (or typically with finite correlation-length), we recall that the following optimal estimates of \eqref{eq:optima-sp-exp} hold (see \cite{Gloria-Neukamm-Otto-a}): for all $\nu>0$,
\begin{equation*}
\int_0^\nu  \d e_{\mfk{d}}(\lambda) \lesssim \nu^{d/2+1}.
\end{equation*}
In view of \eqref{2-1:eq:KV-2}, such a control of the behavior of the spectral measure at the edge of the spectrum should give information on the speed of convergence of $t^{-1} \E_0\Ll[R_t^2\Rr]$ to $0$. Indeed, one can show using basic calculus manipulations (see the proof of \cite[Proposition~8.2]{Mourrat-10}) that \eqref{eq:optima-sp-exp} implies that in fact, $\E_0\Ll[R_t^2\Rr]$ remains bounded as $t$ tends to infinity (or grows no faster than $\log(t)$ in dimension $2$). The crux of the proof is thus to show that \eqref{eq:optima-sp-exp}. 
In the case of discrete time, our estimates of the spectral are slightly weaker (we have a logarithmic divergence in dimension 2, and optimal estimates up to dimension 6), see \cite[Appendix~A]{Gloria-Mourrat-10a}. This is however sufficient to prove Theorem~\ref{t:MCdisc}.

In the decomposition \eqref{decompsquare}, there finally remains to study the cross-product $\E_0\Ll[ M_t R_t \Rr]$. Since we know that roughly speaking, $\E_0\Ll[R_t^2\Rr]$ remains bounded as $t$ tends to infinity, a naive Cauchy-Schwarz estimate would ensure us that $|\E_0\Ll[ M_t R_t \Rr]|$ grows no faster than $\sqrt{t}$. But in fact, a second ``miracle'' occurs here, since this cross product turns out to be $0$ ! Recall that it can be written as
\begin{equation}
\label{vanish}
\E_0\Ll[ \Ll(\xi \cdot X_t + \phi(X_t,\omega) - \phi(0,\omega) \Rr)\Ll(\phi(X_t,\omega) - \phi(0,\omega)\Rr) \Rr].
\end{equation}
To see that this expectation is $0$, it suffices to show that its derivative at time $0$ vanishes, thanks to the Markov property and to the stationarity of $(\omega(t))_{t \ge 0}$. Using \eqref{nonconstr}, we obtain that this derivative is equal to
$$
\expec{  \sum_{z \sim 0} \omega_{0,z}\Ll(\xi \cdot z + \phi(z,\omega) - \phi(0,\omega)\Rr) \Ll(\phi(z,\omega) - \phi(0,\omega) \Rr)  } = \expec{A (\xi + \nabla \phi) \cdot \nabla \phi },
$$
and the latter is equal to $0$ by \eqref{2-1:eq:corr-sto} (the fact that \eqref{vanish} vanishes is also clear from the alternative construction of the corrector based on orthogonal projections, as was done for instance in \cite{matpiat}).

This completes our sketch of the proof of the continuous-time analog of Theorem~\ref{t:MCdisc}. For discrete time, we refer the reader to \cite{Gloria-Mourrat-10b}.

\subsubsection{Methods based on other functions of the random walk}

We now turn our attention to numerical methods based on the convergence in \eqref{convergegendisc}, and in fact, we continue to discuss instead the continuous-time analog \eqref{convergegen}. We focus on a special case of \eqref{convergegen}, namely, the fact that for any bounded continuous $f : \R^d \to \R$, 
\begin{equation}
\label{converget}
\E_0\Ll[ f\Ll(\frac{X_t}{\sqrt{t}}\Rr) \Rr] \xrightarrow[t \to + \infty]{} E[f(B_1)],
\end{equation}
where $B_1$ is a Gaussian vector with covariance matrix $2\Ah$.
In the preceeding chapter, we highlighted two ``miracles'' that made the difference between the l.h.s.\ of \eqref{convergesquare} and its limit as small as $t^{-1}$. Both these miracles are very special of the square function. Although we cannot prove it, for a generic functional, we expect the error to be of order at least $t^{-1/2}$ in general. We now present the known upper bounds on the error.

\begin{thm}
\label{t:quantclt}
Let $B_1$ be a Gaussian vector with covariance matrix $2 \Ah$.

(1) \cite{Mourrat-12} If $f$ is smooth and bounded, then for any $\delta > 0$, there exist $q,C$ such that
\begin{equation}
\label{e:quantclt}
\Ll|\E_0\Ll[ f\Ll(\frac{X_t}{\sqrt{t}}\Rr) \Rr] - E[f(B_1)] \Rr| \le C
\left|
\begin{array}{ll}
	t^{-1/4} & \text{if } d = 1, \\
	\log^q(t) \ t^{-1/4} & \text{if } d = 2, \\
	t^{-1/2+\delta} & \text{if } d \ge 3.
\end{array}
\right.
\end{equation}

(2) \cite{Mourrat-11a} There exists $q,C$ such that for any $\xi \in \R^d$,
\begin{equation}
\label{e:quantclt2}
\sup_{x \in \R}\Ll| \P_0\Ll[ \frac{\xi \cdot X_t}{\sqrt{t}} \le x   \Rr]  - P\Ll[ \frac{\xi \cdot B_1}{\sqrt{t}} \le x   \Rr]\Rr| \le C 
\left|
\begin{array}{ll}
	t^{-1/10} & \text{if } d = 1, \\
	\log^q(t) \ t^{-1/10} & \text{if } d = 2, \\
	\log(t) \ t^{-1/5} & \text{if } d = 3, \\
	t^{-1/5} & \text{if } d \ge 4.
\end{array}
\right.
\end{equation}
\end{thm}
Part (2) of the theorem corresponds to considering the function $f(z) = \1_{\xi\cdot z \le x}$, and is often known as a \emph{Berry-Esseen} estimate.

For comparison, note that if $X_t$ was replaced by a sum of i.i.d.\ centred random variables with finite third moment and covariance matrix $2\Ah$, then the left-hand sides of \eqref{e:quantclt} and \eqref{e:quantclt2} would both be $O(t^{-1/2})$ (see for instance \cite[Theorem~5.5]{Petrov}). 
It was recently realized (but yet to be written) that the bound in \eqref{e:quantclt} for $d = 2$ can be improved to $t^{-1/2+\delta}$ for any $\delta > 0$. The conjectured optimal rates in this inequality are
\begin{equation}
\label{conjrates}
\Ll|
\begin{array}{ll}
t^{-1/4} & \text{if } d = 1, \\
\log(t) \ t^{-1/2} & \text{if } d = 2, \\
t^{-1/2} & \text{if } d \ge 3.	
\end{array}
\Rr.
\end{equation}

\begin{rem}
\label{rem:iid}
To be precise, if $X_t$ was replaced by a sum of i.i.d.\ centred random variables, say $Z_1, \ldots, Z_t$ ($t \in \N$ here), then the characteristic function of the normalized sum would satisfy (for any $\zeta \in \R$)
$$
\E\Ll[ \exp\Ll(i \zeta \frac{Z_1 + \cdots + Z_t}{\sqrt{t}}\Rr)  \Rr] = \E\Ll[ \exp\Ll(i\zeta  \frac{Z_1}{\sqrt{t}}\Rr)\Rr]^{t} = \exp\Ll( -\frac{\zeta^2\E Z_1^2}{2}  - i \frac{\zeta^3 \E Z_1^3}{\sqrt{t}} + O(t^{-1})  \Rr),
$$	
where we assumed that $Z_1$ has finite moments of order $4$. The first term in the expansion gives the characteristic function of the limiting Gaussian random variable, and there is indeed a correction of order $t^{-1/2}$. But \emph{this correction vanishes} in the special case when $\E Z_1^3 = 0$, and in particular if the distribution of $Z_1$ is invariant under the transformation $z \mapsto -z$. Numerical simulations suggest that similar cancellations sometimes also happen for the random walk in random environment, as we will discuss below.
\end{rem}

Without entering into details, we now sketch the main steps of the proof of Theorem~\ref{t:quantclt}. Given that, roughly speaking, $\E_0[R_t^2]$ remains bounded, it is not too hard to control the error due to the absence of the ``miracle of orthogonality''. We thus turn our attention to the absence of the first mentionned ``miracle'', which consists in the fact that the relation \eqref{miracle1} is exact, regardless of whether the martingale $M_t$ is actually closely resembling a Brownian motion or not. For a generic function $f$, it is no longer true that the difference
\begin{equation}
\label{diff}
\Ll|\E_0\Ll[ f\Ll(\frac{M_t}{\sqrt{t}}\Rr) \Rr] - E[f(B_1)] \Rr|
\end{equation}
vanishes, and the proof of Theorem~\ref{t:quantclt} requires that we find bounds on this quantity. Recall that we derived the fact that the difference in \eqref{diff} converges to $0$ from the convergence of $t^{-1} V_t$ in \eqref{ergodicthm}. The two main steps of the proof of Theorem~\ref{t:quantclt} are:

(i) to estimate the speed at which the variance of $t^{-1} V_t$ tends to $0$ (i.e.\ to get an $L^2$ estimate on the speed of convergence in \eqref{ergodicthm}), and

(ii) to rely on a general quantitative version of the central limit theorem for martingales which transports the estimate obtained in (i) into an upper bound on \eqref{diff}.


At least for $d \ge 4$, our control of the variance of $t^{-1} V_t$ is optimal, and the reason for the not-so-intuitive exponent $1/5$ appearing in part (2) of Theorem~\ref{t:quantclt} is hidden in part~(ii) of the proof, that is, in the general quantitative central limit theorem for martingales. Surprisingly, this general result is however optimal \cite{Mourrat-11b}. Yet, we conjecture that the exponents obtained in part (2) of Theorem~\ref{t:quantclt} are not optimal, and that the rates obtained in part (1) of the theorem may in fact hold as well in part (2).

\subsection{Numerical study}
For our numerical tests, we continue with the case introduced in paragraph~\ref{ss:num1} of i.i.d.\ conductances $(\omega_e)_{e \in \mathbb{B}}$ with 
$$
\P\Ll[\omega_e = 1\Rr] = \P\Ll[\omega_e = 4\Rr] = 1/2.
$$
Recall that by symmetry considerations, we know that in this case, $\Ahd$ is a multiple of the identity. For the two-dimensional case, we know moreover that $\Ah = 2\ \mathrm{Id}$, and thus by \eqref{defAhd}, that $\Ahd = 1/5\  \mathrm{Id}$. We do not know of such a formula when $d = 3$, and for the purpose of the analysis of the methods, we may use an approximation of the homogenized coefficients obtained with the periodization method of paragraph~\ref{ss:period}.

For the simulation of the random walk, we generate the environment only along the trajectory of the walk. More precisely, when the walk arrives at some point, it first checks if the neighboring conductances have already been ``seen''. If so, then the value already generated is used, else a value is taken at random and stored for future usage. Roughly speaking, the walk up to time $n$ discovers of order $n$ conductances, independently of the dimension $d \ge 2$ (except for a logarithmic correction in dimension $2$). 

The methods based on simulating random walks have two main interesting features. One is that their effectiveness is fairly insensitive to the dimension. The other is that the computation is massively parallel, since it requires to simulate a lot of \emph{independent} random walks, each in its own \emph{independent} environment.

\subsubsection{Method based on the mean square displacement}
We start by investigating numerically the method based on the computation of the mean square displacement of the random walk.
As for the methods based on the corrector, we investigate separately the statistical and systematic errors, that is, the two terms in the sum \eqref{decomperrorin2}. For the statistical error, we focus on 
$$
\E_0^\otimes\Ll[\Ll(\hat{A}_k(n) - n^{-1} \td{\E}_0\Ll[ \Ll(\xi \cdot Y_n\Rr)^2\Rr] \Rr)^2\Rr] = \Var{\hat{A}_k(n)},
$$
where we write $\Var{\cdot}$ for the variance with respect to the measure $\P^\otimes_0$. Since $\hat{A}_k(n)$ is a sum of independent random variables, we have
\begin{equation}
\label{scalingvar}
\Var{\hat{A}_{k}(n)} = \frac{1}{k} \Var{\frac{p(\omega)}{\expec{p}} \Ll( \frac{\xi \cdot Y_n}{\sqrt{n}} \Rr)^2}.
\end{equation}
A simple variation of the proof of the convergence in \eqref{convgauss} shows that
$$
\lim_{n \to + \infty} \Var{\frac{p(\omega)}{\expec{p}} \Ll( \frac{\xi \cdot Y_n}{\sqrt{n}} \Rr)^2} = \mathsf{v},
$$
where $\mathsf{v}$ was defined in \eqref{defmathsfv}.
Table~\ref{table:random-err} shows our numerical estimations in dimension $2$ for 
\begin{equation}
\label{e:variance}
\Var{\frac{p(\omega)}{\expec{p}} \Ll( \frac{\xi \cdot Y_n}{\sqrt{n}} \Rr)^2}
\end{equation}
for several values of $n$, and compares it to the predicted limiting value. As we see on \eqref{scalingvar}, the variance of our estimator $\hat{A}_k(n)$ is then obtained by dividing this value by the number $k$ of walks we run. Table~\ref{table:random-err3} presents the same results for the 3-dimensional case.

\begin{table}
\begin{tabular}{|c|c|c|c|c|c|c|c|c|c|}
\hline
$n$ & 10 & 20 & 40 & 80 & 160 & 320 & 640 & 1280 & $\infty$
\\
\hline
Variance & 0.40 & 0.39 & 0.380 & 0.373 & 0.369 & 0.367 & 0.3653 & 0.3647 & 0.3632
\\
\hline
\end{tabular}
\bigskip
\caption{Numerical estimates for \eqref{e:variance} and the theoretical limiting value in dimension $2$.}
\label{table:random-err}
\end{table}

\begin{table}
\begin{tabular}{|c|c|c|c|c|c|c|c|c|c|}
\hline
$n$ & 10 & 20 & 40 & 80 & 160 & 320 & 640 & 1280 & $\infty$
\\
\hline
Variance & 0.20 & 0.19 & 0.188 & 0.186 & 0.184 & 0.184 & 0.1836 & 0.1835 & 0.1829
\\
\hline
\end{tabular}
\bigskip
\caption{Numerical estimates for \eqref{e:variance} and the theoretical limiting value (computed with a numerical estimate for $\Ahd$) in dimension $3$.}
\label{table:random-err3}
\end{table}


%
%
%
%

Turning to the systematic error, the theoretical prediction is that displayed in \eqref{systerr}. To test this, we computed $\hat{A}_{k_n}(n)$ with $k_n = K(n) n^2$, where $K(n)$ is some large number. This choice of $k_n$ ensures that the random fluctuations are washed out, so that $\hat{A}_{k_n}(n)$ is very close to its expectation. In practice, we chose $K(n) = 10^4$ for $n \le 320$, and $K(n) = 10^3$ for larger values. 

\begin{figure}
\centering
\psfrag{hh}{$\rm{log}_{10}(N)$}
\psfrag{kk}{$\rm{log}_{10}$(Syst. error)}
\includegraphics[scale=0.7]{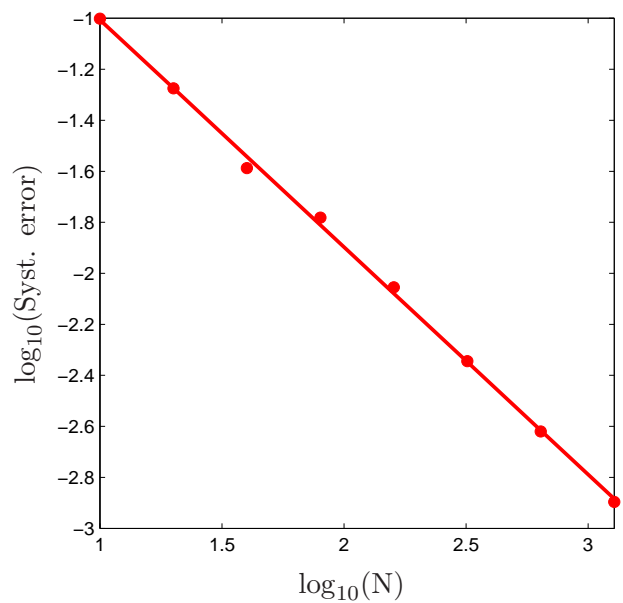}
\caption{The systematic error in two dimensions, rate $0.89$ and prefactor $0.77$.}
\label{fig:syst-err}
\end{figure}

For the two-dimensional case, our results are reported on Figure~\ref{fig:syst-err}. In this case, the theoretical prediction in \eqref{systerr} contains a polylogarithmic correction to the rate of decay in $1/n$. We conjecture that the true rate is actually $\ln(n)/n$. On the log-log plot, this would create a correction to the expected slope of $-1$ of the order of $1/\ln(n)$, which is between $0.17$ and $0.14$ for $n$ between $320$ and $1280$, roughly corresponding to the observed slope.

\begin{figure}
\centering
\psfrag{hh}{$\rm{log}_{10}(N)$}
\psfrag{kk}{$\rm{log}_{10}$(Syst. error)}
\includegraphics[scale=.7]{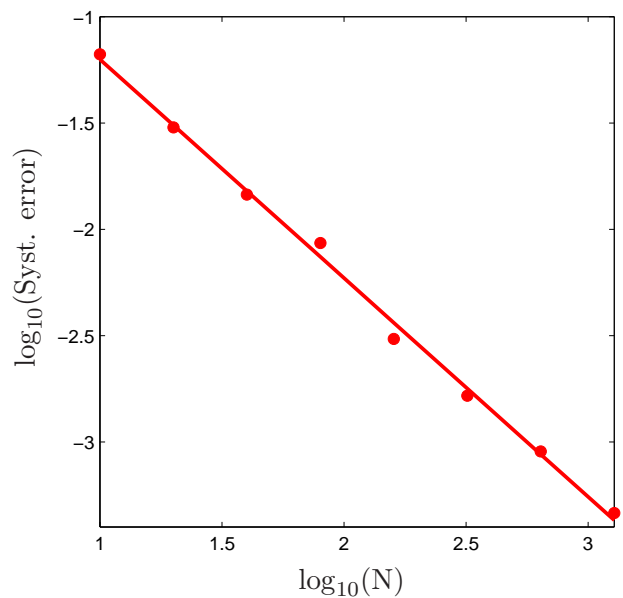}
\caption{The systematic error in three dimensions, rate $1.03$ and prefactor $0.67$.}
\label{fig:syst-err-3d}
\end{figure}

In three dimensions, our findings are reported on Figure~\ref{fig:syst-err-3d}. Here, the observations match very well with our predicted convergence rate of $1$.

\subsubsection{Methods based on other functions of the random walk}
We now turn our attention to methods based on computing other functions of the final position of the random walk than the squared displacement. For any reasonable function $f$, we can devise an estimator as we did for the special function $f(x) = (\xi \cdot x)^2$, which takes the form
\begin{equation}
\label{defhatAf}
\hat{A}^f_k(n) := \frac{1}{k \expec{p}}\sum_{i = 1}^k p(\omega^{(i)}) f(Y^{(i)}_n/\sqrt{n}).
\end{equation}
We have 
$$
\E_0^\otimes[\hat{A}^f_k(n)] = \td{\E}_0[f(Y_n/\sqrt{n})],
$$
which (for suitable $f$) converges to $E[f(B_1)]$, where $B_1$ is a Gaussian vector with covariance matrix $2\Ahd$. 
Since this more general estimator is still a sum of i.i.d.\ random variables, the statistical error is easy to understand, as we have 
$$
\Var{ \hat{A}^f_k(n) } = \frac{1}{k} \Var{ \frac{p(\omega)}{\expec{p}} f\Ll( \frac{\xi \cdot Y_n}{\sqrt{n}} \Rr) },
$$
where the last variance converges (for suitable $f$) to a constant that can be explicitely written.
We now focus our numerical investigation to the systematic error, that is, the difference between $\td{\E}_0[f(Y_n/\sqrt{n})]$ and the limiting value $E[f(B_1)]$. Based on the analogy with the continuous-time case, we expect that the systematic error will decay as $n^{-1/2}$ for $d \ge 2$, possibly with a logarithmic correction in dimension $2$. 

In order to test this prediction, we chose the function
\begin{equation}
\label{deff}
f(x) = \exp\Ll( -\frac{\|x\|_2^2}{2} \Rr).
\end{equation}
In the case we investigate numerically, the conductances being i.i.d.,\ we know that the homogenized matrix is of product form, say $2\Ahd = \sigma^2 \ \mathrm{Id}$, and actually we also know that $\sigma^2 = 2/5$. A simple computation shows that in this case and for $f$ defined in \eqref{deff}, the limiting value is given by
$$
E[f(B_1)] = (\sigma^2+1)^{-d/2},
$$
which is equal to $5/7$ in dimension $2$, and which we can approximate thanks to our previous numerical estimation of $\sigma^2$ in dimension $3$.

\begin{figure}
\centering
\psfrag{hh}{$\rm{log}_{10}(N)$}
\psfrag{kk}{$\rm{log}_{10}$(Syst. error)}
\includegraphics[scale=.7]{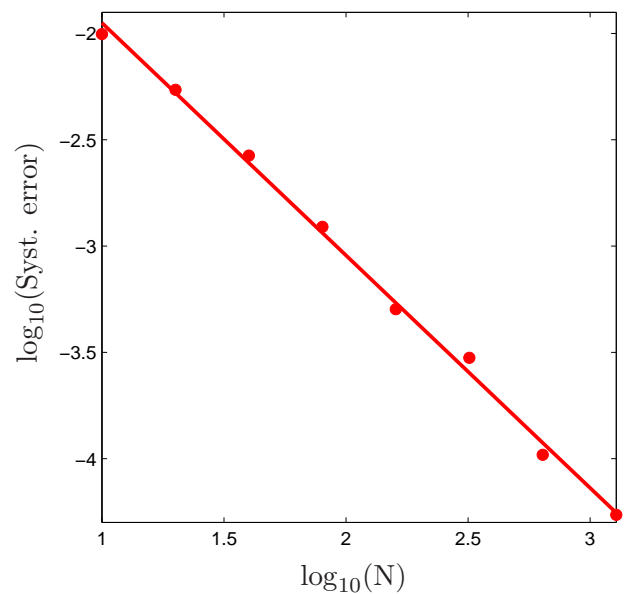}
\caption{The systematic error in two dimensions for $f$ as in \eqref{deff}, rate $1.09$ and prefactor $0.14$.}
\label{fig:sys-err-f-2}
\end{figure}

\begin{figure}
\centering
\psfrag{hh}{$\rm{log}_{10}(N)$}
\psfrag{kk}{$\rm{log}_{10}$(Syst. error)}
\includegraphics[scale=.7]{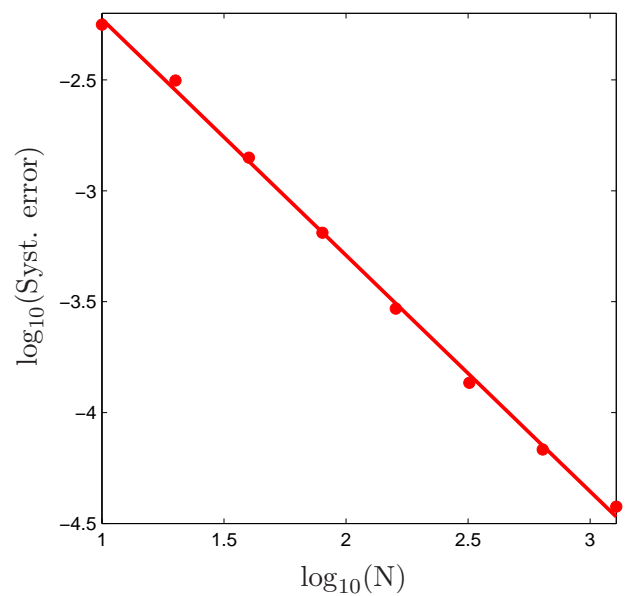}
\caption{The systematic error in three dimensions for $f$ as in \eqref{deff}, rate $1.06$ and prefactor $0.07$.}
\label{fig:sys-err-f-3}
\end{figure}

The systematic errors we obtained numerically are shown on Figures~\ref{fig:sys-err-f-2} and \ref{fig:sys-err-f-3} for dimensions $2$ and $3$ respectively. 

Surprisingly, the observed convergence rates are far better than the predicted ones, being close to $1$ in both cases instead of the predicted $1/2$. Similar rates where observed for other choices of the function $f$. 

We believe that these surprising rates are due to special cancellations which do not hold for arbitrary distributions. As was explained already in Remark~\ref{rem:iid}, such cancellations happen already in certain cases when one considers sums of i.i.d.\ random variables, in particular when the distribution of the random variables is invariant under the transformation $z \mapsto -z$. Another way to guess what the correct order of the correction is to compute an asymptotic expansion in terms of the ellipticity ratio, in the spirit of what was done in \cite[Appendix]{Gloria-Otto-09}. If one looks at the first order correction, it turns out that it behaves as $t^{-1/2}$ in general, but cancellations occur if the distribution is symmetric under the transformation $z \mapsto -z$, and the correction becomes of order $t^{-1}$. In view of these two facts, we believe that the convergence rates obesrved on Figures~\ref{fig:sys-err-f-2} and \ref{fig:sys-err-f-3} are somehow miraculous, and due to the symmetry of the distribution of the conductances under the transformation $z \mapsto -z$.

\begin{figure}
\centering
\includegraphics[scale=.8]{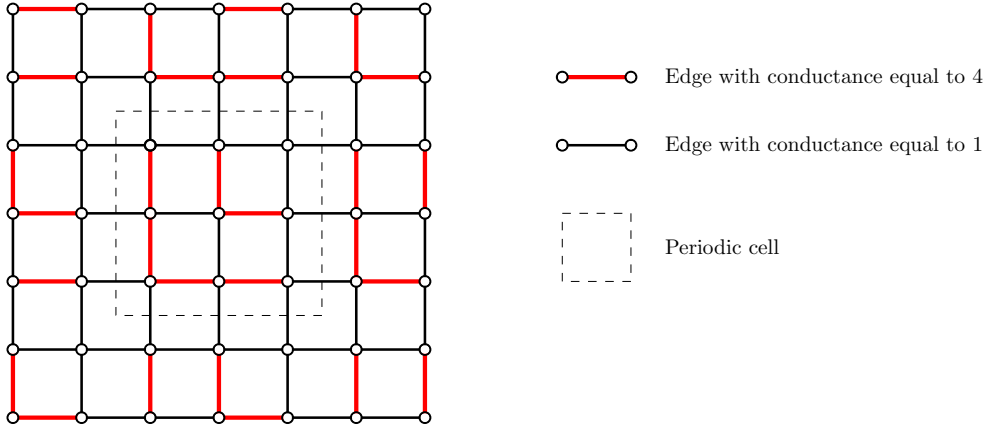}
\caption{The periodic environment, without $z \mapsto -z$ symmetry}
\label{fig:cell}
\end{figure}

\begin{figure}
\centering
\psfrag{hh}{$\rm{log}_{10}(N)$}
\psfrag{kk}{$\rm{log}_{10}$(Syst. error)}
\includegraphics[scale=.7]{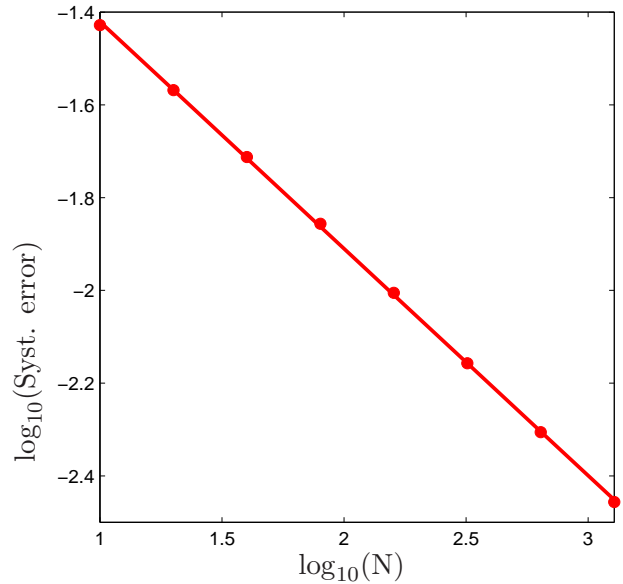}
\caption{The systematic error in the periodic environment for $f(x,y) = \sin(x)$, rate $0.49$ and prefactor $0.12$.}
\label{fig:sys-err-perio}
\end{figure}

In order to convince ourselves that the convergence rates observed are indeed non-typical, we reverted to a study of a periodic environment, choosing it so that it is not symmetric under the transformation $z \mapsto -z$. A periodic environment should provide better convergence rates than any generic and truly random environment (in any dimension). We chose a simple $3$-periodic cell displayed on Figure~\ref{fig:cell}. We chose $f(x,y) = \sin(x)$, so that $\td{\E}_0[f(Y_n/\sqrt{n})]$ tends to $0$ as $n$ tends to infinity (and we thus do not need to compute the homogenized matrix to study the systematic error). The results we obtained are shown on Figure~\ref{fig:sys-err-perio}. As expected, they display a convergence rate of $1/2$.

Finally, we investigated the convergence rates for the non-smooth function 
\begin{equation}
\label{deffsup}
f(x) = \1_{\xi \cdot x \le z},
\end{equation}
where we took $\xi$ to be the first vector of the canonical basis and $z = 1/2$ or $z = 1/4$. For the systematic error, Theorem~\ref{t:quantclt} ensures a convergence rate of $1/10$ in dimension $2$, and of $1/5$ in dimension $3$, up to logarithmic corrections. We believe that the exponent can be pushed to $1/5$ in dimension $2$ with a refined argument (yet to be written), but the proof cannot be pushed to higher exponents \cite{Mourrat-11b}. 

\begin{figure}
\centering
\psfrag{hh}{$\rm{log}_{10}(N)$}
\psfrag{kk}{$\rm{log}_{10}$(Syst. error)}
\includegraphics[scale=.7]{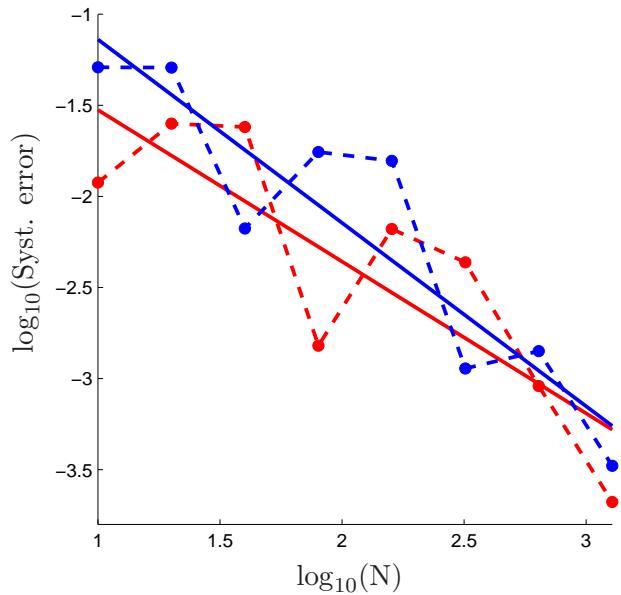}
\caption{The systematic error in two dimensions for $f$ as in \eqref{deffsup} and $z = 1/2$ (red), rate $0.8$ and prefactor $0.2$, or $z = 1/4$ (blue), rate $1.0$ and prefactor $0.7$.}
\label{fig:sys-err-2D-sup}
\end{figure}

\begin{figure}
\centering
\psfrag{hh}{$\rm{log}_{10}(N)$}
\psfrag{kk}{$\rm{log}_{10}$(Syst. error)}
\includegraphics[scale=.7]{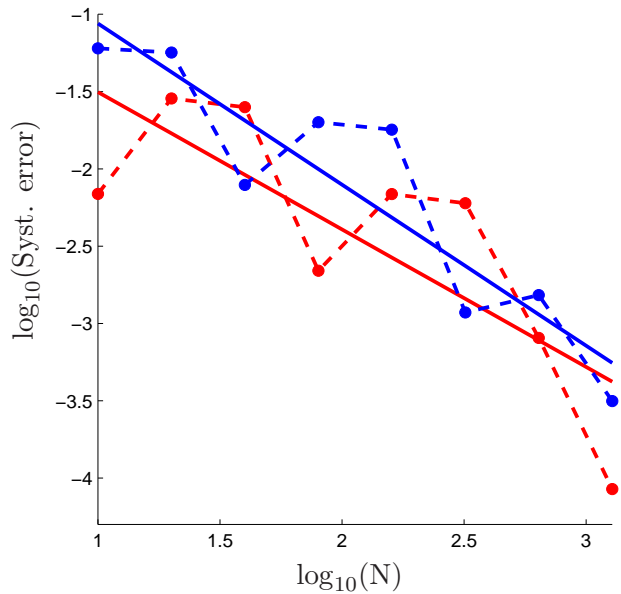}
\caption{The systematic error in three dimensions for $f$ as in \eqref{deffsup} and $z = 1/2$ (red), rate $0.9$ and prefactor $0.2$, or $z = 1/4$ (blue), rate $1.0$ and prefactor $1.0$.}
\label{fig:sys-err-3D-sup}
\end{figure}

The main goal of these numerical investigations is to see whether this limitation to the exponent $1/5$ is an artefact of the method of proof. Our results are displayed on Figure~\ref{fig:sys-err-2D-sup} for $d=2$, and on Figure~\ref{fig:sys-err-3D-sup} for $d = 3$. 

While the onset of a nice asymptotic regime is delayed, the results suggest that indeed the convergence rates ultimately settle to a behavior similar to that observed on Figure~\ref{fig:sys-err-f-2} or \ref{fig:sys-err-f-3}, that is, ultimately displaying a convergence rate close to the value $1$.

\section{Conclusion}

In this article we have recalled qualitative and quantitative results on stochastic homogenization of discrete linear elliptic PDEs and of random walks in random environments on $\Z^d$.
This has allowed us to make a rather complete picture of numerical methods to approximate homogenized coefficients, based both on the corrector equation and on the random walk.
Numerical tests have confirmed the sharpness of the analysis, supported some conjectures, and put some interesting phenomena into evidence (such as ungenerically fast decay rates due to specific symmetries of the environement).
We hope this contribution will help mathematicians identify challenging conjectures, and help practitioners make mathematically-based choices on the method to use in more concrete cases.

\medskip

We have only considered discrete elliptic equations.
The discrete case has the advantage of being numerically inexpensive to simulate compared to the case of continuum linear elliptic equations (for which we have to appeal to approximation methods such as the finite element, finite difference or finite volume methods).
The extension of the results of this paper to the continuum case is currently under investigation, and 
\cite{Gloria-Otto-09,Gloria-Otto-09b} have already been extented to the continuum case \cite{Gloria-Otto-10b}. 
The extension of the RWRE approach to the continuum space is yet to be done.
We hope to address these issues in future works.

\end{document}